\input amstex
\documentstyle{amsppt}
\pageno=1
\magnification1200
\catcode`\@=11
\def\logo@{}
\catcode`\@=\active
\NoBlackBoxes
\vsize=23.5truecm
\hsize=16.5truecm

\def\d{d\!@!@!@!@!@!{}^{@!@!\text{\rm--}}\!}

\def\crp{\overline{\Bbb R}_+}

\def\rnp{{\Bbb R}^n_+}

\def\crnp{\overline{\Bbb R}^n_+}

\def\ang#1{\langle {#1} \rangle}

\def\ttilde{\overset{\,\approx}\to}

\def\Zfrac{\tsize\frac1{\raise 1pt\hbox{$\scriptstyle z$}}}

\def\rp{ \Bbb R_+}

\define\tr{\operatorname{tr}}

\define\Tr{\operatorname{Tr}}

\def\res{\operatorname {res}}
\def\Res{\operatorname {Res}}
\def\ord{\operatorname {ord}}

\topmatter
\title
Trace expansions and the noncommutative residue for manifolds with
boundary
\endtitle
\author{Gerd Grubb and Elmar Schrohe}\endauthor
\address Department of Mathematics, University of Copenhagen,
Universitetsparken 5, \linebreak DK-2100 Copenhagen, Denmark.\endaddress 
\email grubb\@math.ku.dk\endemail

\address Universit\"at Potsdam, Institut f\"ur Mathematik, Postfach
60 15 53, D-14415 Potsdam, Germany.\endaddress 
\email schrohe\@math.uni-potsdam.de\endemail
\abstract For a pseudodifferential boundary operator $A$ of order $\nu\in\Bbb Z $
and class 0 (in the Boutet de Monvel calculus) on a compact
$n$-dimensional manifold with boundary, we consider the function
$\operatorname{Tr}(AB^{-s})$, where $B$ is an auxiliary system
formed of the Dirichlet realization of a second order strongly elliptic
differential operator and an elliptic operator on the boundary. We
prove that 
$\operatorname{Tr}(AB^{-s})$ has a meromorphic
extension to $\Bbb C$ with poles at the half-integers $s= (n+\nu-j
)/2$, $j\in\Bbb N$
(possibly double for $s<0$), and we prove that
its residue at 0 equals
the noncommutative residue of $A$, as defined by Fedosov, Golse,
Leichtnam and Schrohe by a different method. To achieve this, we
establish a full asymptotic expansion of
$\operatorname{Tr}(A(B-\lambda )^{-k})$ in powers $\lambda ^{-l/2}$ and
log-powers $\lambda ^{-l/2}\log\lambda $, where the noncommutative
residue equals the coefficient of the highest order log-power. There
is a related expansion of $\operatorname{Tr}(Ae^{-tB})$.

The paper will appear in Journal Reine Angew\. Math\. (Crelle's Journal).
\endabstract
\rightheadtext{Noncommutative residue}
\endtopmatter

\document
\subhead 1. Introduction\endsubhead

In this paper we study trace expansions associated to operators in
the calculus of Boutet de Monvel on a compact manifold $X$ with boundary. 
As an auxiliary operator, we first fix an invertible system $B$ formed by the
Dirichlet realization of a strongly elliptic differential operator in the
interior of $X$ and an elliptic operator on the boundary. 
Given a pseudodifferential boundary value
problem $A$  in the Boutet de Monvel calculus, the mapping  
$$s\mapsto \Tr(AB^{-s}),$$
then is a holomorphic function of $s$ for large $\operatorname{Re}s$. 
We show that it extends to a
meromorphic function on the whole complex plane with at most double poles.
Moreover, we prove that the noncommutative residue
$\res(A)$ of $A$ can be recovered as a residue in this expansion: 
$$
\res(A) = \ord B\cdot \Res_{s=0}\Tr (AB^{-s}),\tag{1.1}$$
provided that $A$ is of class zero. In particular, the right hand
side does not 
depend on the choice of $B$.

In order to put these results
into perspective, let us recall the situation of a
compact manifold {\it without} boundary. In 1984, Wodzicki found a trace on the
algebra $\Psi(X)$ of all classical pseudodifferential operators on a closed
manifold $X$ of dimension $n>1$. In fact, this traces lives on the
quotient $\Psi(X)/\Psi^{-\infty}(X)$ modulo the regularizing elements, and it
is the unique trace there --- up to multiples --- if $X$ is
connected.  (A {\it 
trace} on an algebra ${\Cal A}$ is a linear map $\tau: {\Cal A}\to {\Bbb C}$
such that $\tau(ab)= \tau(ba)$ for all $a,b\in {\Cal A}.$)
Wodzicki called it the {\it noncommutative residue}.

His approach was based on an analogue of  (1.1): 
Given an arbitrary classical pseudodifferential operator $A$ acting on sections of a vector bundle $E$
over $X$, he chose an invertible classical pseudodifferential operator $P$ of sufficiently large positive
order (larger than ${\ord}\,A$), 
satisfying Agmon's condition for the existence of a ray of minimal
growth for the resolvent. 
Following Seeley [S67], this allows the construction of the complex powers 
$(P+uA)^s$, $s\in{\Bbb C}$, for $u\in{\Bbb R}$ close to zero, and a meromorphic extension of ${\Tr}(P+uA)^{-s}$
to ${\Bbb C}$. 
Wodzicki proved that the right hand side of the formula, below, is independent
of $P$ and defined
$$
\res(A) = \ord P\cdot\left.
\frac d{du}\right|_{u=0} \Res_{s=1}\Tr\left((P+uA)^{-s}\right).
\tag{1.2}
$$

It is clear from Seeley's results that the above
expression only depends on finitely many
terms in the asymptotic expansion $\sum a_j(x,\xi)$ of
the symbol  $A$ into terms $a_j$ which are
homogeneouos of degree $j$ in $\xi$.
Wodzicki showed that in fact 
$$
\res(A) = (2\pi)^{-n} \int_{S^*X}\tr_E a_{-n}(x,\xi)\,d\sigma;
\tag{1.3}
$$
here $\tr_E$ is a trace in $\operatorname{Hom}(E)$ and $d\sigma$ 
the surface measure on $S^*X$.  

The notion `noncommutative residue' is justified also by a historical reason. 
In
the late seventies Manin [M79] and Adler [A79] in their work on the
Korteweg-de Vries
equation introduced a trace functional on an algebra of formal pseudodifferential 
operators
in one dimension.
The symbols were Laurent series in one variable with coefficients in $\Bbb C$,
and to the operator 
with the symbol $\sum a_j\xi^j$ they associated
the trace $a_{-1}$, i.e., the usual 
(`commutative') residue of complex analysis. 
As formula (1.3) shows, 
Wodzicki's residue is a natural generalization
of the functional employed by Manin and Adler to the $n$-dimensional case. 
    
The noncommutative residue was discovered independently by Guillemin [Gu85] 
in connection with his alternative proof of Weyl's law on the asymptotic 
distribution of the eigenvalues of an elliptic operator. An interesting survey
of these facts was given by Kassel in [K89].

For a compact manifold $X$  with boundary $\partial X=X'$, 
the situation is more complicated.
As pointed out by Wodzicki, there
is no trace $(\not=0)$ on the algebra of all classical pseudodifferential operators with the Leibniz product 
on $X$ unless $\partial X=\emptyset $: 
This algebra is not sensitive to what happens at the boundary.
A more natural choice in this context is the Boutet de Monvel algebra [BM71] 
of polyhomogeneous 
pseudodifferential boundary value problems 
(or pseudodifferential boundary  operators, 
$\psi$dbo's for short) of integer order. 
An operator in this algebra is a map $A$ acting on sections of vector bundles 
$E$ over $X$ 
and $F$ over $X'$, given by a matrix 
$$
A=\pmatrix P_++G & K\\ \quad\\ T&S\endpmatrix\: 
\matrix
C^\infty (X,E)& & C^\infty (X,E)\\
\oplus&\to&\oplus\\
C^\infty (X',F)& & C^\infty (X',F)
\endmatrix
\tag{1.4}
$$
(Here $P_+$ is a $\psi $do truncated to $X$, $G$ is a singular Green
operator (s.g.o\.), $K$ a Poisson operator, $T$ a trace operator and
$S$ a $\psi $do on $X'$.)
Indeed it could be shown by Fedosov, Golse, Leichtnam, and Schrohe
[FGLS96] that there is a 
 trace on this algebra, also called res. It is given by
$$
\multline
\operatorname{res}(A)
={(2\pi )^{-n}}\int_{S^*X}\tr_E p_{-n}(x,\xi )\,d\sigma  \\
+{(2\pi )^{1-n}}
\int_{S^*X'}[\tr_E (\tr_n g)_{1-n}(x',\xi ')+
\tr_F s_{1-n}(x',\xi ')]\,d\sigma '
,
\endmultline
\tag1.5
$$
(we have corrected a sign error in [FGLS96, (2.19)]).
Here $d\sigma'$ is the surface measure on $S^*X'$ and $\tr_n\,g$ stands for the {\it normal trace} 
of the symbol $g$ of the singular Green operator  $G$: 
$$
(\tr_n g)(x',\xi')=(2\pi)^{-1}\int^+g(x',\xi',\xi_n,\xi_n)\,d\xi_n,
\tag{1.6}
$$
it is a polyhomogeneous pseudodifferential symbol on $X'$. Assuming
that $X$ is 
connected, `res' turns out to be the unique continuous trace
 --- up to multiples --- on the 
quotient modulo the regularizing operators. 

In [FGLS96],
the definition of $\res (A)$ and the proof of the trace property
relied completely on 
the symbolic level, since the analytic properties of traces of powers of
$\psi$dbo's were not known in sufficient generality. It remained an open
problem to identify the noncommutative residue with a residue in the sense of
(1.1) or (1.2). 

Progress was possible as a result of the work of Grubb and Seeley [GS95] on 
weakly parametric pseudodifferential operators on closed manifolds. 
In fact, [GS95, Theorem 2.7] shows that, for classical pseudodifferential
operators $A$ and $P$ of integer order,
where $P$ is elliptic with parameter in a sector in $\Bbb C$, there
are expansions (regardless of the magnitude of $\ord P$):

$$\alignedat2 
\text{(I)}&&\quad \Tr\bigl(A(P-\lambda
)^{-k}\bigr)&\sim \sum_{j\ge 0}c_j(-\lambda )
^{\frac{n+\operatorname{ord}A-j}{\operatorname{ord}P}-k}+\sum_{l\ge
0}(c'_l\log (-\lambda )+c''_l)(-\lambda )
^{-l-k},\\
\text{(II)}&&\quad \Tr\bigl(Ae^{-tP}\bigr)&\sim \sum_{j\ge 0}\tilde
c_j t
^{\frac{j-n-\operatorname{ord}A}{\operatorname{ord}P}}
+\sum_{l\ge 0}(-\tilde c'_l\log t +\tilde c''_l)t
^{l} ,\\
\text{(III)}&&\quad \Gamma (s)\Tr\bigl(AP^{-s})&\sim \sum_{j\ge 0}\frac{\tilde
c_j }
{s+\frac{j-n-\operatorname{ord}A}{\operatorname{ord}P}}
+\sum_{l\ge 0}\Bigl(\frac{\tilde c'_l}{(s+l)^2} +\frac{\tilde
c''_l}
{s+l}\Bigr)
 .
\endalignedat\tag1.7$$
In (I), $k>(n+\operatorname{ord}A)/\operatorname{ord}P$, and
$\lambda $ goes to infinity on suitable rays in $\Bbb C\setminus \Bbb R_+$; 
in (II), $t\to 0+$, and it holds when $P$ is strongly elliptic; 
in (III), the sign ``$\sim$'' indicates that the left hand side is
meromorphic on $\Bbb C$ with pole structure as in the right hand side. 
The coefficients $\tilde c_j, \tilde c'_l,\tilde c''_l$ are
proportional to the coefficients $c_j,  c'_l, c''_l$ by
universal constants. 
In particular, $$
\operatorname{Res}_{s=0}\Tr(AP^{-s})=\tilde c'_0=c'_0.$$
We recall moreover that the statements (I), (II)
and (III) are essentially equivalent, as 
explained e.g\. in [GS96], via transition formulas such as
$$\gathered
P^{-s}=\tfrac{1}{(s-1)\cdots(s-k)}\tfrac{ i}{2\pi
}\int_{{\Cal C}}\lambda 
^{k-s}\partial _\lambda ^k(P-\lambda )^{-1} \,d\lambda=
\tfrac1{\Gamma (s)}\int_0^\infty t^{s-1}e^{-tP}\,dt,\\ 
e^{-tP}=t^{-k}\tfrac{ i}{2\pi }\int_{{\Cal C'}}e^{-t\lambda }
\partial _\lambda ^k(P-\lambda )^{-1} 
\,d\lambda=\tfrac1{2\pi i}\int_{\operatorname{Re}s=c}t^{-s}\Gamma
(s)P^{-s}\,ds;
\endgathered\tag 1.8$$
note that $\partial _\lambda
^{k}(P-\lambda )^{-1}=k!(P-\lambda )^{-k-1}$.

The point of view of determining the residue from formulas (1.5),
together with more general expansion
theorems for
parameter-dependent $\psi $do's in [GS95], is 
developed by
Lesch in [L99], who extends the considerations to a graded algebra
of logarithmically polyhomogeneous $\psi $do's, defining higher
residues as well.

In the present paper we show how the results of [GS95] can be applied to the
analysis of boundary value problems in the Boutet de Monvel calculus. 
It is the main technical difficulty to establish expansions analogous to 
(I), (II), and (III) in (1.7) above. 
In order to describe our result, let $A$ be
an arbitrary
$\psi$dbo of order $\nu\in\Bbb Z$ and class zero as in (1.4); we
denote $\operatorname{dim}E=n'$.
We next choose our auxiliary operator. We let $P_1$ be a strongly
elliptic second order differential operator in $E$, which has scalar
principal symbol on a neighborhood of $X'$; by $P_{1,D}$ we denote its
Dirichlet realization. Finally we pick a strongly elliptic second order
pseudodifferential operator $S_1$ on $C^\infty(X',F)$ 
and set 
$
B=\pmatrix P_{1,D}& 0\\ 0&S_1\endpmatrix.
$
Then 
$$
\Tr (A(B-\lambda)^{-k})
= \Tr_X \left((P_++G)(P_{1,D}-\lambda)^{-k}\right)
+\Tr_{X'}\left(S(S_1-\lambda)^{-k}\right)\tag1.9
$$
is well-defined for $k>(n+\nu)/2$. The behavior of 
$\Tr_{X'}\left(S(S_1-\lambda)^{-k}\right)$ is covered by the preceding 
results: We have an expansion as in (I) of (1.7) with $n$ replaced
by $n-1$; moreover, the coefficient of the first logarithmic term in this
expansion is
proportional to the noncommutative residue of $S$ on $X'$. 
It therefore remains to analyze the first term. This the content of the 
following theorem:

\proclaim{Theorem 1.1} 
There are full asymptotic expansions (when $k>(n+\nu )/2$)
$$\alignedat2
\text{\rm (I)}&&\; \Tr((P_++G)(P_{\operatorname{1,D}}-\lambda
)^{-k})&\sim \sum_{j\ge 0}\!c_{j}(-\lambda)
^{\frac{n+\nu - j}{2}-k}+\sum_{l\ge
0}(c'_l\log (-\lambda )+c''_l)(-\lambda )
^{-\frac l2-k},\\
\text{\rm (II)}&&\;
\Tr((P_++G)e^{-tP_{1,\operatorname{D}}})&\sim
\sum_{j\ge 0}\tilde  
c_{j} t
^{\frac{j-n-\nu  }{2}}
+\sum_{l\ge 0}(-\tilde c'_l\log t +\tilde c''_l)t
^{\frac l2},\\
\text{\rm (III)}&&\; \Gamma (s)\Tr
((P_++G)P_{1,\operatorname{D}}^{-s})&\sim \sum_{j\ge 0}\frac{\tilde
c_j } 
{s+\frac{j-n-\nu  }{2}}
+\sum_{l\ge 0}\Bigl(\frac{\tilde c'_l}{(s+\frac l2)^2} +\frac{\tilde
c''_l}
{s+\frac l2}\Bigr)
 .
\endalignedat
\tag1.10$$
The coefficients $\tilde c_j, \tilde c_l', \tilde c''_l$ are
proportional to the coefficients $c_j,  c_l', c''_l$ by
universal constants.
In particular,
$\operatorname{Res}_{s=0}\Tr((P_++G)P_{1,\operatorname{D}}^{-s})$ equals
$\tilde c'_0=c'_0$, and $$
\operatorname{res}(P_++G)=\operatorname{ord}P_1\cdot
\operatorname{Res}_{s=0}\Tr((P_++G)P_{1,\operatorname{D}}^{-s}).\tag1.11
$$
\endproclaim 

The main effort in our proof lies in the deduction of (I) in (1.10).
Besides this, we show how the specific contributions to (1.5)
arise in the trace expansion of
$(P_++G)(P_{1,\operatorname{D}}-\lambda )^{-k}$.

\example{Remark 1.2}
Our method also applies to the case where 
$G$ is of class $r>0$. In general, the expansion (I) will then take
the form
$$
\Tr((P_++G)(P_{\operatorname{1,D}}-\lambda
)^{-k})\sim \sum_{j\ge 0}\!c_{j}(-\lambda)
^{\frac{n+\nu - l}{2}-k}+\sum_{l\ge
-r}(c'_l\log (-\lambda )+c''_l)(-\lambda )
^{-\frac l2-k},\tag 1.12
$$
and the other expansions are accordingly modified.
The relation (1.11) will not in general be
valid. For example, if $G=K\gamma _0$, res$(G)$ is generally nonzero
according to [FGLS96], whereas $G(P_{1,D}-\lambda )^{-k}$ is zero since
$(P_{1,D}-\lambda )^{-1}$ maps into a space where $\gamma _0u$
vanishes. (We denote $u|_{X'}=\gamma _0u$.)
\endexample

Let us moreover remark that the Dirichlet condition may be replaced
by the Neumann condition or other coercive boundary conditions, since the
resulting singular Green operators have symbol structures similar to
the case studied in detail here.

Partial expansions were obtained  for 
large classes of parameter-dependent operators in [G96] and its
predecessor from 1986 (see also the
appendix of [G92]), however without getting as far as the
residue at zero or the logarithmic
term that is interesting in the present context. Full expansions
were given for some problems connected with Dirac operators in
[GS95], [G99], where the 
efforts were concentrated on reducing the question to a study of $\psi
$do's on the boundary, to which the calculus
of so-called weakly polyhomogeneous $\psi $do's could be applied.
An extension of the calculus to a class of parameter-dependent $\psi$dbo's 
has been worked out in [G00] for an algebra of operators containing the
strongly polyhomogeneous $\psi$dbo's. It provides full trace expansions,
but it does not include the general compositions with parameter-independent
operators that we need to treat here.  

An interesting point in the present investigation is that the 
important contribution comes from a logarithmic term, whereas the
logarithmic terms are considered more as  
a disturbance in the problems connected with Dirac operators.
\smallskip

{\it Plan of the paper:} In Section 2 we write
$(P_++G)(P_{1,D}-\lambda )^{-k}$ as a sum of five terms of different
nature, reduce the problem to local coordinates at the boundary, and
determine the structure of the symbols connected with the auxiliary
operator as rational functions. Section 3 begins with a simple
example demonstrating the basic idea of the proof; then we use Laguerre
expansions to extend this idea to the two terms composed of singular
Green operators. Sections 4 and 5 treat the terms composed of a $\psi
$do and an s.g.o. In Section 6, the various results are collected to
give the  proof of Theorem 1.1. The special Laguerre expansions we
use are recalled in the Appendix.

\subhead 2. The structure of the compositions \endsubhead

We are going to establish an asymptotic expansion for $\Tr(AR_\mu ^k)$,
where $A=P_++G$ and $R_\mu =(P _{1,\operatorname{D}}+\mu^2)^{-1}$. 

Our main tool is the Boutet de Monvel calculus. 
We assume familiarity with the standard notions in this context, 
recommending the reference [G96].
Under reasonable conditions on 
$P_1$, one can show that $R_\mu$, acting on $C^\infty(X,E)$,
is a strongly parameter-dependent 
polyhomogeneous element of the Boutet de Monvel calculus,
of the kind of the upper left corner of the matrix in (1.4). 
Its pseudodifferential part simply is $((P_1+\mu^2)^{-1})_+$, 
so that the usual parametrix construction gives us full information 
on its symbol.
The analysis of the singular Green part of $R_\mu$ 
requires considerably more care. 
Knowing the structure of its symbol components
will enable us, however, to also determine the
symbolic structure of the $k$-th power $R_\mu^k$ and, eventually, that
of $AR^k_\mu$.  

In order to show the trace expansions we then rely on the concepts
of polyhomogeneous parameter-dependent symbol classes; references here 
are [G96] and [GS95]. Let us just emphasize one aspect:
When symbols depend on a parameter 
$\mu $ (running in
a sector $\Gamma $ of $\Bbb C\setminus\{0\}$), 
we distinguish between a weak and a strong version of
the homogeneity property $p(x,t\xi ,t\mu )=t^m p(x,\xi ,\mu )$. When
it holds for all $\xi ,\mu ,t$ with $|\xi |^2+|\mu |^2\ge 1$, $t\ge
1$, $(\xi ,\mu )\in \Bbb R^n\times (\Gamma \cup\{0\})$, $p$ is said to
be {\it strongly homogeneous} (of degree $m$).
When
it holds for all $\xi ,\mu ,t$ with $|\xi |\ge 1$, $t\ge
1$, $(\xi ,\mu )\in \Bbb R^n\times \Gamma $, $p$ is said to
be {\it weakly homogeneous} or just homogeneous. In the weakly
homogeneous case, one must further specify which behavior is assumed
in the noncompact set $\{\,(\xi ,\mu )\mid |\xi |\le 1,\,|\mu |\ge 1\,\}$.
Strong resp\. weak {\it poly}-homogeneity of a symbol means that it has an
asymptotic expansion 
in strongly resp\. weakly homogeneous terms of decreasing orders,
with specific requirements on the remainder estimates. Briefly
explained,
{\it strong
polyhomogeneity} means that the new variable $\mu $ enters on a par with
the cotangent variables $\xi $, so that both $\mu $-derivatives 
and $\xi _j$-derivatives give decrease in the estimates in terms of
$\ang{(\xi ,\mu )}$. (For $y=(y_1,\dots ,y_l)$, we denote $(|y_1|^2+\dots
+|y_l|^2+1)^{\frac12}=
\ang{y}$.) In the trace expansions obtained below, strongly
polyhomogeneous operators contribute no logarithmic terms.

Let us now describe the set-up more precisely: 
Without loss of generality we may suppose that $P_1$ is the
restriction of a strongly elliptic second-order
differential operator $P_1$ defined in a bundle $\widetilde E$
extending $E$ over a closed manifold $\widetilde
X$ neighboring $X$, with scalar principal symbol on a neighborhood of $X'$.
$P_+$ can 
be assumed to be the truncation to $X$ of a $\psi $do $P$ of order
$\nu $ defined on 
$\widetilde X$.
We fix $k>(n+\nu)/2$ so that $AR^k_\mu$ is trace class.
We also assume that the lower bounds of  $P_1$ (on $\widetilde X$) as well as 
$P_{1,\operatorname{D}}$ (on $X$) are positive. 
As explained in [G96, Section 3.3], 
there exists an obtuse keyhole region$$
W=\{z\in\Bbb C \mid|z |\le r\text{ or }|\arg z|\le\tfrac\pi
2+\varepsilon\}\tag2.1
$$
(for some $\varepsilon >0$) such that when $\mu^2 \in W$, 
the inverses $Q_\mu =(P_1+\mu ^2)^{-1}$
(on $\widetilde X$) and $R_\mu =(P _{1,\operatorname{D}}+\mu
^2)^{-1}$ (on $X$) exist as bounded operators in $L_2$ and are $O(\mu
^{-2})$ on rays in $W$.
Moreover, by [G96, (3.3.21)], 
$$
R_\mu =Q_{\mu ,+}+G_\mu ,\tag2.2
$$ 
where $Q_{\mu ,+}=(Q_\mu)_+$ is the restriction of $Q_\mu$
to $X$ and $G_\mu $ is a singular Green operator, described in
detail in Lemma 2.4 below. 
Then $$\aligned
R_\mu ^k&=(Q_{\mu ,+}+G_\mu )^k=(Q_\mu ^k)_{+}+G^{(k)}_\mu ,\text{
with}\\
G^{(k)}_\mu &=(Q_{\mu ,+})^k-(Q_\mu ^k)_++\operatorname{pol}(G_\mu
,Q_{\mu ,+});\endaligned\tag2.3 
$$here $\operatorname{pol}(G_\mu ,Q_{\mu ,+})$ is a
linear combination of compositions of the (non-commuting) factors $G_\mu $ and $Q_{\mu ,+}$
with $k$ factors in each term, at least one of them being $G_\mu $.
The symbol structure of $G^{(k)}_\mu$ will be given in Proposition 2.5. 
So, writing $(Q^k_\mu )_+=Q^k_{\mu ,+}$, $$
AR_\mu ^k=(P_++G)(Q^k_{\mu ,+}+G^{(k)}_\mu )=(PQ^k_\mu )_+-L(P,Q^k_\mu
)+P_+G^{(k)}_{\mu } 
+GQ^k_{\mu ,+}+GG^{(k)}_\mu .\tag 2.4
$$
Here, $L(P,Q^k_\mu)= (PQ^k_\mu)_+-P_+Q^k_{\mu,+}$ is the leftover term in the composition
of $P_+$ and $Q^k_{\mu,+}$. All these operators are trace class since $k>(n+\nu)/2$. 
There are five terms to be dealt with:

$$
\Tr((PQ^k_\mu )_+),\quad -\Tr(L(P,Q^k_\mu )),\quad
\Tr(P_+G^{(k)}_\mu ),
\quad \Tr(GQ^k_{\mu ,+}),\quad
\Tr(GG^{(k)}_\mu );\tag 2.5
$$
they will now be analyzed in detail.

The first term in (2.5) is covered by [GS95, Th\. 2.1], applied to
the extensions of 
$P$ and $Q^k_\mu $ to $\widetilde X$, when we integrate the resulting 
kernel formula over $X$.
It remains to consider the other four. It will be relatively easy to deal with the effects 
caused away from the boundary:

\subsubhead{2.a Contributions from the interior} \endsubsubhead

We shall first show that the only interesting contributions stem from a neighborhood of the 
boundary:
Using a partition of unity $1=\sum_{i=1}^I\theta _i(x)$  subordinate
to a cover of $X$ by coordinate patches $(\Omega _j)_{j=1}^J$
(with trivializations of the bundles) in such
a way that any four of the functions $\theta _i$ are supported in one
of the coordinate patches (see
e.g\. [GK93, Lemma 1.5]), we can localize any of the operators in
(2.5). For those that are products $\Cal A_{1}\Cal A_2$, we
write$$
\Cal A_1\Cal A_2=\sum_{i,j,l,m\le I}\theta _{i}\Cal A_1\theta
_j\theta _l\Cal A_2\theta _m,\tag2.6
$$which reduces the problem to the consideration of
compositions $\Cal A'_{1}\Cal A'_{2}$ where both factors act in one
coordinate patch $\Omega $
 and map into compactly supported functions there. 
We can assume that the principal symbol of $P_1$ is a scalar times
the identity matrix in those patches that meet the boundary.

If $\overline\Omega $ is disjoint from the boundary, the factor 
$\theta _lG_\mu ^{(k)}\theta _m$ coming from the third and the fifth
expression in (2.5) is strongly 
polyhomogeneous and of order $-\infty $. By [G96, Theorem 2.5.11],
a strongly polyhomogeneous operator of order $-r$ is
continuous from $H^{s,\mu }(X,E)$ to $H^{s+r,\mu }(X,E)$ for $s> -\frac12$,
hence its composite with an operator of
order $\nu $ is trace class when $r>n+\nu $ with a trace norm that is  
 $O(\mu ^{n+\nu -r})$ for $\mu \to\infty $ in $\Gamma $,
$$
\Gamma =\{\mu  \in\Bbb C \mid |\arg \mu 
 |\le\tfrac\pi
4+\tfrac12\varepsilon \}\tag 2.7
$$
(cf\. (2.1)). So the terms with the factor $\theta _lG_\mu
^{(k)}\theta _m$  are trace class operators
whose traces are $O(\mu ^{-N})$ for any $N$,
hence do not contribute to the trace expansions. 

For the second expression in (2.5), we can assume that the $\theta _i$ form a
partition of unity on $\widetilde X$, so that we can write
$$\multline
L(P_+,Q^k_\mu )=(PQ^k_\mu )_+-P_+Q^k_{\mu ,+}\\
=\sum_{i,j,l,m\le
I}(\theta _{i}P
\theta
_j\theta _lQ^k_\mu \theta _m)_+ -
\sum_{i,j,l,m\le
I}(\theta _{i}P
\theta
_j)_+(\theta _lQ^k_{\mu} \theta _m)=\sum_{i,j,l,m\le
I}L(\theta _{i}P
\theta
_j,\theta _lQ^k_\mu \theta _m).\endmultline
$$
In the local coordinates for the coordinate patch $\Omega $ where
$\theta _i,\theta _j,\theta _l,\theta _m$ are supported, 
$$
L(\theta _{i}P
\theta
_j,\theta _lQ^k_\mu \theta _m)=G^+(\theta _{i}P
\theta
_j)G^-(\theta _lQ^k_\mu \theta _m),
$$
by [G96, (2.6.27)]. Here, when  $\overline\Omega $ is disjoint from
$X'$, $G^-(\theta _lQ^k_\mu \theta _m)$ is a strongly polyhomogeneous 
s.g.o\. of order $-\infty $, so the product is trace class with $O(\mu ^{-N})$ estimates
for any $N$.

In the fourth expression in (2.5), when $\overline\Omega $ is
disjoint from $X'$,  
the operator $\theta _iG\theta _j$ is a singular Green operator of
order $-\infty $, hence has a $C^\infty $ kernel, and identifies with
a negligible $\psi $do $\Cal P_{ij,+}$. Here
$$
\theta _iG\theta _j\theta _lQ^k_{\mu ,+}\theta _m=
\Cal P_{ij,+}\theta _lQ^k_{\mu ,+}\theta _m=
(\Cal P_{ij}\theta _lQ^k_{\mu }\theta _m)_+-
L(\Cal P_{ij},\theta _lQ^k_{\mu }\theta _m).
$$
The $L$-term gives $O(\mu ^{-N})$ estimates just like the
preceding case. For $\Cal P_{ij}\theta _lQ^k_{\mu }\theta _m$
we have from the calculus of [GS95] for $\psi $do's on $\widetilde X$
that it belongs to $S^{-\infty ,-2k}$ there, hence has a diagonal
kernel expansion $$
\sum_{r\in \Bbb N}c_r(x)\mu ^{-2k-r}\text{ for $\mu \to\infty $ in
}\Gamma ,\tag2.8$$ 
which, when
integrated over $X$, gives the trace expansion $\sum_{r\in \Bbb N}c_r\mu
^{-2k-r}$.

So we see that all the interior contributions from the second to
fifth term in (2.5) have
trace expansions that are either void or have the form (2.8),
with no logarithms.  

We can therefore focus on a neighborhood of the boundary where 
we can work in local coordinates on $\crnp$. 
We shall use the notation in (2.5), omitting explicit mention of the
$\theta _i$'s or coordinate changes. 

\subsubhead{2.b The symbol of $Q_\mu^k$} \endsubsubhead

Let  us first recall the symbol structure of $Q_\mu^k$.
For a boundary patch where $P_1$
has principal symbol $p_{1,2}I$ with a scalar $p_{1,2}$, the
symbol of $Q_{\mu }$ has the 
following form in local coordinates:
$$\aligned
q(x,\xi ,\mu )&\sim \sum_{l\in\Bbb N}q_{-2-l}(x,\xi ,\mu ),\text{ with}\\
q_{-2}(x,\xi ,\mu )&=(p_{1,2}(x,\xi )+\mu ^2)^{-1}I,\\
q_{-2-J}(x,\xi ,\mu )&=\sum_{2\le m\le 2J+1}\frac{r_{J,m}(x,\xi
)}{(p_{1,2}(x,\xi )+\mu ^2)^{m}},\text{ for }J>0;\endaligned\tag 2.9
$$
here the $r_{J,m}$ are $n'\times n'$-matrices ($n'=\operatorname {dim} E$) of
homogeneous polynomials in $\xi $ of degree $2m-J$ with smooth
coefficients  (cf\. [S67, (1)--(1a)]). 

For $Q_\mu ^k$ we find from the composition rules for $\psi $do's that
the symbol $q^{\circ k}$ satisfies:$$\aligned
q^{\circ k}(x,\xi ,\mu )&\sim \sum_{l\in\Bbb N}q^{\circ k}_{-2k-l}(x,\xi ,\mu ),\text{ with}\\
q^{\circ k}_{-2k}(x,\xi ,\mu )&=(p_{1,2}(x,\xi )+\mu ^2)^{-k}I,\\
q^{\circ k}_{-2k-J}(x,\xi ,\mu )&=\sum_{k+1\le m\le 2J+k}\frac{r_{k,J,m}(x,\xi
)}{(p_{1,2}(x,\xi )+\mu ^2)^{m}},\text{ for }J>0;\endaligned\tag 2.10
$$the $r_{k,J,m}$ being matrices of
homogeneous polynomials in $\xi $ of degree $2m-2k-J$. 
\pagebreak

\subsubhead{2.c The symbols of $G^\pm(Q_\mu^k)$} \endsubsubhead

The second term in (2.5) is $L(P,Q^k_\mu )$ which can be written 
$L(P,Q^k_\mu )=G^+(P)G^-(Q^k_\mu )$, cf. [G96, (2.6.25--26)].
We shall now study the structure of $G^-(Q^k_\mu )$. 
For simplicity of notation we shall write $x'$ instead of $(x',0)$.

The polynomial 
$$p_{1,2}(x',\xi )+\mu ^2=a(x')\xi _n^2+b(x',\xi ')\xi
_n+c(x',\xi ')+\mu 
^2\tag 2.11$$
 has two roots with respect to $\xi _n$,
$$
i\kappa ^+(x',\xi ',\mu
)\in\Bbb C_+\quad\text{ and 
}\quad-i\kappa ^-(x',\xi ',\mu )\in\Bbb C_-,\tag 2.12
$$ strongly homogeneous of degree
$1$ with
$\operatorname{Re}\kappa ^{\pm}>0$ when $\mu ^2\in W$, cf\.
(2.1), so
$$\aligned
p_{1,2}(x',\xi )+\mu ^2&=a(x')(\xi _n-i\kappa ^+(x',\xi ',\mu ))
(\xi _n+i\kappa ^-(x',\xi ',\mu ))\\
&=
a(x')(\kappa ^+(x',\xi ',\mu )+i\xi _n)
(\kappa ^-(x',\xi ',\mu )-i\xi _n).\endaligned\tag 2.13
$$
Then we have formulas (by decomposition of rational functions of $\xi
_n$ in simple fractions)
$$\multline
\frac1{(p_{1,2}(x',\xi )+\mu ^2)^{m}}=\frac1{a(x')^{m}(\kappa
^+(x',\xi ',\mu )+i\xi
_n)^{m}(\kappa ^-(x',\xi ',\mu )-i\xi
_n)^{m}}\\
=a(x')^{-m}\Bigl(\sum_{1\le j\le m}\frac{a^+_{m,j}(x',\xi ',\mu
)}{(\kappa ^+(x',\xi ',\mu )+i\xi 
_n)^{j}}+\sum_{1\le j\le m}\frac{a^-_{m,j}(x',\xi ',\mu )}{(\kappa ^-(x',\xi ',\mu )-i\xi
_n)^{j}}\Bigr),\endmultline\tag 2.14
$$
where the $a_{m,j}^\pm(x',\xi ',\mu )$ are strongly homogeneous of degree
$j-2m$ in $(\xi ',\mu )$.

Recall that when a rational function $f(\xi _n)$ with
poles in $\Bbb C\setminus\Bbb R$ is decomposed as $f=f^++f^-+f'$,
where $f^{\pm}$ have poles in $\Bbb C_{\pm}=\{\xi _n\in\Bbb C\mid
\operatorname{Im}z\gtrless 0\}$ and are $O(\ang{\xi _n}^{-1})$, and $f'$
is polynomial, then the projections $h^+$, $h^-$, $h^-_{-1}$ (cf\. [G96,
Prop\. 2.2.2]), are simply the mappings $h^+\:f\mapsto
f^+$, $h^-\:f\mapsto f^-+f'$, $h^-_{-1}\:f\mapsto f^-$. (2.14) gives for $m=k$
the decomposition $q^{\circ k}_{-2k}=h^+q^{\circ k}_{-2k}+ h^-q^{\circ
k}_{-2k}$; here $h^-q^{\circ k}_{-2k}=h^-_{-1}q^{\circ k}_{-2k}$.

We can also apply $h^+$ and $h^-$ to the lower order terms in
$q^{\circ k}$. Since
$q^{\circ k}_{-2k-J}(x',\xi ,\mu )$ is a proper rational function of
$\xi _n$ with poles $i\kappa ^+$ and $-i\kappa ^-$ of order $\le 2J+k$
for $J\ge 0$ (cf\. (2.10)), and all coefficients as well as $\kappa ^+$
and $\kappa 
^{-}$ are strongly polyhomogeneous in $(\xi ',\mu )$, one finds by
decomposition in simple fractions that $q^{\circ k}_{-2k-J}(x',\xi
,\mu )$ is the sum of a proper rational function $h^+q^{\circ
k}_{-2k-J}(x',\xi ,\mu )=q^{\circ k,+}_{-2k-J}(x',\xi ,\mu )$ of $\xi
_n$ with poles $i\kappa
^+$ and a proper rational function $h^-q^{\circ
k}_{-2k-J}(x',\xi ,\mu )=q^{\circ k,-}_{-2k-J}(x',\xi ,\mu )$ of $\xi
_n$ with
poles $-i\kappa ^-$ (as described in detail in (2.15) below); the
coefficients are strongly homogeneous in $(\xi ',\mu )$ and the poles are of orders $\le 2J+k$. The functions
$q^{\circ k\pm}_{-2k-J}$ are homogeneous of degree $-2k-J$ in $(\xi
,\mu )$. Thus we find:

\proclaim{Proposition 2.1} For $x=(x',0)$, denoted $x'$,
there is a decomposition of $q^{\circ k}$ into a part
$q^{\circ k,+}=h^+q^{\circ k}$ with poles in
$\Bbb C_+$ and a part $q^{\circ k,-}=h^-q^{\circ k}$ with poles in $\Bbb C_-$
(as functions of $\xi _n$) of the
following form:
$$\aligned
\quad q^{\circ k}=q^{\circ k,+}+q^{\circ k,-},&\quad
q^{\circ k,\pm}(x',\xi ,\mu )\sim\sum_{
J\ge 0}q^{\circ k,\pm}_{-2k-J}(x',\xi ,\mu ),\\
q^{\circ k}_{-2k-J}(x',\xi ,\mu )&=q^{\circ k,+}_{-2k-J}(x',\xi ,\mu
)+q^{\circ k,-}_{-2k-J}(x',\xi ,\mu ),\\
q^{\circ k,\pm}_{-2k-J}(x',\xi ,\mu )&=\sum_{1\le j\le
2J+k}\frac{r^{\pm}_{k,J,j}(x',\xi ',\mu )}
{(\kappa ^\pm(x',\xi ',\mu )\pm i\xi
_n)^{j}},
\endaligned\tag 2.15
$$where the $n'\times n'$-matrix formed
numerators $r^{\pm}_{k,J,j}(x',\xi ',\mu )$ are strongly
homogeneous in $(\xi ',\mu )$ of degree $j-2k-J$. 

Moreover, in the Taylor expansion of $q^{\circ k}$,$$
\sum_{l\in\Bbb N}\tfrac{x_n^l}{l!}\partial _{x_n}^lq^{\circ
k}(x',0,\xi ,\mu ),\tag 2.16
$$
the terms satisfy:
$$\aligned
\quad \partial _{x_n}^lq^{\circ k}=\partial _{x_n}^lq^{\circ
k,+}+\partial _{x_n}^lq^{\circ k,-},&\quad 
\partial _{x_n}^lq^{\circ k,\pm}(x',\xi ,\mu )\sim\sum_{
J\ge 0}\partial _{x_n}^lq^{\circ k,\pm}_{-2k-J}(x',\xi ,\mu ),\\
\partial _{x_n}^lq^{\circ k}_{-2k-J}(x',\xi ,\mu )&=\partial _{x_n}^lq^{\circ k,+}_{-2k-J}(x',\xi ,\mu
)+\partial _{x_n}^lq^{\circ k,-}_{-2k-J}(x',\xi ,\mu ),\\
\partial _{x_n}^lq^{\circ k,\pm}_{-2k-J}(x',\xi ,\mu )&=\sum_{1\le j\le
2J+k+l}\frac{r^{l,\pm}_{k,J,j}(x',\xi ',\mu )}
{(\kappa ^\pm(x',\xi ',\mu )\pm i\xi
_n)^{j}},
\endaligned\tag 2.17
$$where
the numerators $r^{l,\pm}_{k,J,j}(x',\xi ',\mu )$ are strongly
homogeneous of degree $j-2k-J$. 
\endproclaim

The indication ``$\circ k$'' will be left out when $k=1$.
 
\demo{Proof} (2.15) has already been shown. For (2.17) we note as a
starting point that (2.10) implies$$\multline
\partial _{x_n}q^{\circ k}_{-2k-J}(x,\xi ,\mu )
=\sum_{k+1\le m\le
2J+k}\bigl(\frac{\partial _{x_n}r_{k,J,m}}{(p_{1,2}+\mu ^2)^{m}}-
\frac{r_{k,J,m}m\partial _{x_n}p_{1,2}}{(p_{1,2}+\mu ^2)^{m+1}}\bigr) 
\\=\sum_{k+1\le m\le
2J+k+1}\frac{r^{1}_{k,J,m}(x,\xi
)}{(p_{1,2}(x,\xi )+\mu ^2)^{m}},
\endmultline
$$
and hence, successively$$
\partial _{x_n}^lq^{\circ k}_{-2k-J}(x,\xi ,\mu )
=\sum_{k+1\le m\le
2J+k+l}\frac{r^{l}_{k,J,m}(x,\xi
)}{(p_{1,2}(x,\xi )+\mu ^2)^{m}},
$$
for all $l$, with homogeneous polynomials $r^{l}_{k,J,m}(x,\xi )$
of degree $2m-2k-J$ in $\xi $. At $x_n=0$ we  decompose these expansions just
as above, which leads to formulas as in (2.15), except that
the index $j$ runs up to $2J+k+l$.\qed
\enddemo

The symbol-kernel and symbol of
$G^-(Q^k_{\mu })$ are found from the $h^-$-projection of the symbol of
$Q^k_\mu $ by the formulas in [G96, Th\. 2.6.10 and 2.7.4].  
(We recall the relation between an s.g.o\. symbol $g(\xi _n,\eta _n)$ of class 0 and the associated symbol-kernel $\tilde
g(x_n,y_n)$:$$
\tilde g(x_n,y_n )=r^+_{x_n} r^+_{y_n}\Cal F^{-1}_{\xi
_n\to x_n
}\overline{\Cal F}^{-1}_{\eta _n\to y_n}
g(\xi 
_n,\eta _n ),\; g(\xi 
_n,\eta _n )
=\Cal F_{x_n\to\xi _n }\overline{\Cal F}_{y_n\to\eta
_n }e^+_{x_n} e^+_{y_n}\tilde g(x_n,y_n);
\tag 2.18
$$
here  $r^+$  denotes restriction from
$\Bbb R$ to $\Bbb R_+$, and $e^+$ denotes extension 
by 0.) 
As a variant of (A.7), $\overline{\Cal F}^{-1}_{\xi _n\to
z_n}(\kappa ^--i\xi
_n)^{-j}=H(z_n)\frac1{(j-1)!}z_n^{j-1}e^{-z_n\kappa ^-}$. Then 
[G96, (2.6.44--45)] give by direct application to (2.15):

\proclaim{Proposition 2.2} If the coefficients
in $P_1$ are independent of  
$x_n$ near $x_n=0$, then
the singular Green operator $G^-(Q^k_\mu )$ has
symbol-kernel  $\tilde g^{\circ k,-}(x',x_n,y_n,\xi ',\mu )$ (defined for
$x_n$ and $y_n\ge 0$) 
and symbol $g^{\circ k,-}(x',\xi ',\xi _n,\eta _n,\mu ) $ derived from the
symbol of 
$q^{\circ k}$ at $x_n=0$, cf\. {\rm (2.15)}:$$\aligned 
&\tilde g^-(q^{\circ k})=\tilde g^{\circ k,-}(x',x_n,y_n,\xi ',\mu )\sim \sum_{J\in\Bbb
N}\tilde g^{\circ k,-}_{-1-2k-J}(x',x_n,y_n,\xi ',\mu ),\text{ with}\\ 
&\tilde g^{\circ k,-}_{-1-2k-J}(x',x_n,y_n,\xi ',\mu )=\sum_{1\le j\le
2J+k}{r^{-}_{k,J,j}(x',\xi ',\mu )}\tfrac
1{(j-1)!}(x_n+y_n)^{j-1}e^{-(x_n+y_n)\kappa ^-}
,\\
&\qquad=\sum_{1\le j\le
2J+k}{r^{-}_{k,J,j}(x',\xi ',\mu )}\sum_{0\le j'\le
j-1}\tfrac{1}{j'!}x_n^{j'}e^{-x_n\kappa 
^-}\tfrac1{(j-1-j')!}y_n^{j-1-j'}e^{-y_n\kappa ^-}
,\\
&g^-(q^{\circ k})=
g^{\circ k,-}(x',\xi ',\xi _n,\eta _n,\mu )\sim \sum_{J\in\Bbb
N}g^{\circ k,-}_{-1-2k-J}(x',\xi ',\xi _n,\eta _n,\mu ),\text{ with}\\ 
&g^{\circ k,-}_{-1-2k-J}(x',\xi ',\xi _n,\eta _n,\mu )\\
&\qquad=\sum_{1\le j\le
2J+k}{r^{-}_{k,J,j}(x',\xi ',\mu )}
\sum_{0\le j'\le j-1}
\tfrac1{(\kappa ^-(x',\xi ',\mu )+ i\xi
_n)^{j'+1}}
\tfrac1{(\kappa ^-(x',\xi ',\mu )- i\eta 
_n)^{j-j'}};
\endaligned\tag 2.19
$$
since the $r^-_{k,J,j}$ are strongly homogeneous of degree $j-2k-J$ in
$(\xi ',\mu )$, the degree of $g^{\circ k,-}_{-1-2k-J}(x',\xi
',\xi _n,\eta _n,\mu )$ is $-1-2k-J$.

There are similar formulas for the symbol-kernel $\tilde g^+(q^{\circ
k})$ and symbol $ g^+(q^{\circ
k})$ of $G^+(Q^k)$, where $\kappa ^-$ and 
$r^-_{k,J,j}$ are replaced by $\kappa ^+$ and $r^+_{k,J,j}$.
\endproclaim 

If the coefficients in $P_1$ are not independent of $x_n$ near
$x_n=0$, the symbols of 
$G^{\pm}(Q^{\circ k})$ contain moreover the contributions from the terms
in the Taylor expansion (2.16) of $q^{\circ k}$ in $x_n$ at $x_n=0$.

\proclaim{Proposition 2.3}
When the coefficients
in $P_1$ depend on  
$x_n$,
$$
g^\pm(q^{\circ k})(x',\xi,\eta_n,\mu)=g^{\circ k,\pm}(x',\xi ,\eta
_n,\mu )\sim\sum_{j\in\Bbb
N}\tfrac1{j!}\overline D_{\xi_n}^jg^\pm[\partial_{x_n}^jq^{\circ
k}(x',0,\xi,\mu)],\tag2.20 
$$
where $g^\pm$ of the $x_n$-independent symbols are found as in {\rm
(2.19)}, now using the formulas {\rm (2.17)}. 
In particular, in the term of homogeneity degree $-1-2k-J$, the
number of powers
of $\kappa ^{\pm}+ i\xi _n$ and $\kappa ^{\pm}- i\eta _n$ taken
together is at most $2J+k+1$, so the structure is
$$\aligned
 g^{\circ k,\pm}(x',\xi ,\eta _n,\mu )&\sim \sum_{J\in\Bbb
N} g^{\circ k,\pm}_{-1-2k-J}(x',\xi ,\eta _n,\mu ),\text{ with}\\ 
 g^{\circ k,\pm}_{-1-2k-J}(x',\xi ,\eta _n,\mu )&=\sum_{\Sb
j,j'\ge 1,\\j+j'\le
2J+k+1\endSb}\frac{s^{\pm}_{k,J,j,j'}(x',\xi ',\mu )}
{(\kappa ^\pm(x',\xi ',\mu )+ i\xi
_n)^{j}(\kappa ^\pm(x',\xi ',\mu )- i\eta 
_n)^{j'}},
\endaligned\tag2.21$$
with the $s^\pm_{k,J,j,j'}$ strongly homogeneous of degree $j+j'-J-2k-1$.

\endproclaim 

\demo{Proof} We apply [G96, Th\. 2.7.4], which gives the general
formula (2.20). 
Now look at the powers that occur in the resulting 
denominators. As noted in (2.17), $\partial _{x_n}^l$ augments the
range of the 
powers in $q^{\circ k}_{-2k-J}$ to $m\le 2J+k+l$, hence in the
s.g.o\. derived from this, the powers go up to $2J+k+l+1$. The application of
$\overline D_{\xi _n}^l$ furthermore augments the
range of the 
powers to $m\le 2J+k+2l+1=2(J+l)+k+1$. Since the resulting term has
homogeneity degree $-2k-J-l-1$, we find
(2.21) by collecting the terms by degrees of homogeneity.
\qed
\enddemo 

\subsubhead 2.d The symbol of $G_\mu^{(k)}$\endsubsubhead 

As a first step we
shall describe the singular Green operator $G_\mu $ in (2.2)
in a form convenient for the present purposes:

The resolvent $R_\mu $ exists for $\mu ^2\in W$ (e.g\. by variational 
theory), and the solution operator for the fully nonhomogeneous
Dirichlet problem$$
\pmatrix P_1+\mu ^2\\\gamma _0 
\endpmatrix ^{-1}=\pmatrix R_\mu &K_\mu  
\endpmatrix \tag2.22
$$
exists since $\gamma _0 $ has a convenient right inverse. (Recall
that $\gamma _0u=u|_{X'}$.)

We can assume that a normal coordinate $x_n$ and a 
normal derivative 
$D_n$ ($=-{  i}\partial _{x_n}$) have been chosen in a
neighborhood of $\partial X=X'$ so that, with $\varrho =\{\gamma
_0,\gamma _0D_n\}$, we have  
 Green's formula for $P_1$:$$
(P_1u,v)_{X}-(u,P_1^*v)_{X}=(\Cal A\varrho u,\varrho
v)_{X'},\quad \Cal A=\pmatrix \Cal A_{00}(x',D')&{  i}a(x')\\
{  i}a(x')&0 
\endpmatrix ;\tag
2.23 
$$ here $\Cal A_{00}$ is a first order differential operator on $X'$ and
$a(x')$ is the coefficient of $D_{n}^2$ for $P_1$
at $X'$.
Define the nullspace
$$
Z_{\mu ,+}^s=\{ u\in H^s(X)\mid (P_1+\mu ^2)u=0\text{ on }X\}.\tag2.24
$$
Since $P_1+\mu ^2$ is invertible on $\widetilde X$, one can show that
the functions $u\in Z_{\mu ,+}^s$ are uniquely determined by their
Cauchy data $\varrho u$. In fact, the 
Poisson operator$$
K^+_\mu =- r^+ Q_\mu \pmatrix \tilde\gamma _0^*&D_n^*\tilde\gamma
_0^*\endpmatrix\Cal A\:H^{s-\frac12}(X')\times H^{s-\frac32}(X')\to
H^s(X)\tag2.25
$$
acts as an inverse of $\varrho \: Z^s_{\mu ,+}\to  
H^{s-\frac12}(X')\times H^{s-\frac32}(X')$
(cf\. Seeley [S69] or e.g\. [G96, Ex\. 1.3.5]; here $\tilde\gamma _0^*$ is the adjoint of the restriction operator
$\tilde\gamma 
_0$ from functions on
$\widetilde X$ to functions on $X'$). The operator $C^+_\mu =\varrho
K^+_\mu $ is a pseudodifferential projection in
$H^{s-\frac12}(X')\times H^{s-\frac32}(X')$, the {\it Calder\'on
projector}. Clearly, it is strongly polyhomogeneous.

By (2.22) there is also a mapping $K_\mu $ from $H^{s-\frac12}(X')$
to $Z^s_{\mu ,+}$ giving the solution $u=K_\mu \varphi $ of the
semi-homogeneous Dirichlet problem $(P_1+\mu ^2)u=0$, $\gamma
_0u=\varphi $. Composition with $\gamma _0D_n$ leads to the $\psi $do
$P_\mu =\gamma _0D_nK_\mu $, often called the   {\it Dirichlet-to-Neumann
operator}. By [G71, Sect\. 6], the
blocks in $C^+_\mu = 
\pmatrix C^+_{\mu ,00} &C^+_{\mu ,01} \\C^+_{\mu ,10}& C^+_{\mu ,11} 
\endpmatrix $ are elliptic, and $$
P_\mu \sim (C^+_{\mu ,01})^{\circ(-1)}(I- C^+_{\mu ,00})\tag2.26 
$$
where $(C^+_{\mu ,01})^{\circ(-1)}$ is a parametrix of $C^+_{\mu
,01}$; $P_\mu $ is strongly polyhomogeneous.

Now the solution operator $R_\mu \:f\mapsto u$ for the other semihomogeneous
Dirichlet problem $(P_1+\mu ^2)u=f$, $\gamma _0u=0$, is described as
follows:
Let $v=Q_{\mu ,+}f$, then $z=u-v$ solves$$
 (P_1+\mu ^2)z=0,\quad \gamma _0z=-\gamma _0Q_{\mu ,+}f.
$$
Here $\varrho z=\{-\gamma _0Q_{\mu ,+}f,-P_\mu \gamma _0Q_{\mu
,+}f\}$, so $z$ must equal $K^+_\mu $ of this. Thus$$\aligned
u&=v+z=Q_{\mu ,+}f +K^+_\mu \{-\gamma _0Q_{\mu ,+}f,-P_\mu \gamma _0Q_{\mu
,+}f\}
\\
&=Q_{\mu ,+}f +r^+Q_{\mu}\pmatrix \tilde\gamma _0^*&D_n^*\tilde\gamma
_0^*\endpmatrix
\pmatrix\Cal A_{00}(x',D')&{  i}a(x')\\
{  i}a(x')&0  
\endpmatrix \pmatrix I\\P_\mu  
\endpmatrix 
 \gamma _0Q_{\mu ,+}f\\
&=Q_{\mu ,+}f +r^+Q_{\mu} \tilde\gamma _0^*S_{\mu } 
 \gamma _0Q_{\mu ,+}f +r^+Q_{\mu} D_n^*\tilde\gamma _0^*{  i}a 
 \gamma _0Q_{\mu ,+}f,
\endaligned$$
where $S_{\mu }$ is the $\psi $do $\Cal A_{00}+iaP_\mu $ on $X'$.

So if we define
$$\aligned
K_{0,\mu }&=r^+Q_\mu \tilde\gamma _0^*,\quad
K_{1,\mu }=r^+Q_\mu
D_n^*\tilde\gamma _0^*\quad\text{ (Poisson operators)},\\
T_{0,\mu }&=\gamma _0Q_{\mu,+} \quad\text{ (trace operator of order
$-2$ and class 0),}
\endaligned\tag2.27
$$ 
then $G_\mu $ is a sum of two terms, each composed 
of a Poisson operator ($K_{0,\mu }$ resp\. $K_{1,\mu }$), a
$\psi $do on $X'$ ($S_{\mu } $ resp\. 
${  i}a $) and a trace operator $T_{0,\mu }$; all strongly
polyhomogeneous: 
$$
G_\mu = K_{0,\mu }S_{\mu } T_{0,\mu }+K_{1,\mu }{  i}a T_{0,\mu
}.\tag 2.28
$$   

The symbols $k_0(x',\xi ,\mu )$, $k_1(x',\xi ,\mu )$ and
$t_0(x',\xi ,\mu )$
of the Poisson and trace operators $K_{0,\mu }$, $K_{1,\mu }$ and
$T_{0,\mu }$ in local coordinates are easy to find from (2.15)--(2.17)
 (with $k=1$), by use of the rules in [G96] (see in particular the introduction to Section
2.3 there). When the coefficients in $P_1$ are independent of $x_n$
near $x_n=0$, then
$$\aligned
k_0&\sim \sum_{J\in\Bbb N}k_{0,-2-J},\text{ with }
k_{0,-2-J}=\sum_{1\le j\le
2J+1}\tfrac{r^{+}_{1,J,j}}
{(\kappa ^++i\xi
_n)^{j}};\\
k_1&\sim \sum_{J\in\Bbb N}k_{1,-1-J},\text{ with }
k_{1,-1-J}=\sum_{1\le j\le
2J+1}h^+\tfrac{i\xi _nr^{+}_{1,J,j}}
{(\kappa ^++i\xi
_n)^{j}}=\sum_{1\le j\le
2J+1}\tfrac{{r'}^{+}_{1,J,j}}
{(\kappa ^++i\xi
_n)^{j}},\\
&\qquad {r'}^+_{1,J,j}=\kappa ^+r^+_{1,J,j}+r^+_{1,J,j+1};\\
t_0&\sim \sum_{J\in\Bbb N}t_{0,-2-J},\text{ with }
t_{0,-2-J}=\sum_{1\le j\le
2J+1}\tfrac{r^{-}_{1,J,j}}
{(\kappa ^--i\xi
_n)^{j}}.
\endaligned\tag2.29$$
The $r^\pm_{1,J,j}(x',\xi ',\mu )$ are strongly homogeneous of degree
$j-J-2$.  
The ${r'}^+_{1,J,j}(x',\xi ',\mu )$ are strongly homogeneous of degree $j-J-1$,
found from the $r^+_{1,J,j}$ by use of the fact that $i\xi
_n(\kappa ^++i\xi _n)^{-j}
=(\kappa ^++i\xi _n)^{1-j}+\kappa ^+(\kappa ^++i\xi _n)^{-j}$. 

When the coefficients in $P_1$ depend on $x_n$, one can use [G96,
Lemma 2.7.3] as in the proof of Proposition 2.3 above to see that
the formulas for the Poisson and trace symbols have a similar
structure, only with $r^\pm_{1,J,j}$ and ${r'}^+_{1,J,j}$ replaced by other
strongly homogeneous functions of degree $j-J-2$ resp\. $j-J-1$.

In the compositions needed to find the symbol of (2.28), we note that
in the homogeneous term of degree $-3-J$ (order $-2-J$) in for
example $K_{0,\mu }S_{\mu }T_{0,\mu }$, it is terms from 
$k_{0,-2-J'}$ and $t_{0,-2-J''}$ with $J'+J''\le J$ that enter. The
general structure of the symbol of $G_\mu $ is then:
$$\multline
g(x',\xi ',\xi _n,\eta _n,\mu )\sim
\sum_{J\in\Bbb N}g_{-3-J}(x',\xi ',\xi _n,\eta _n,\mu
),\text{ with}\\
g_{-3-J}(x',\xi ',\xi _n,\eta _n,\mu )=\sum_{ \Sb j\ge 1,\,j'\ge 1\\
j+j'\le 2J+2\endSb}\tfrac{s_{J,j,j'}(x',\xi ',\mu )}
{(\kappa ^+(x',\xi ',\mu )+i\xi
_n)^{j}(\kappa ^-(x',\xi ',\mu )-i\xi
_n)^{j'}},\endmultline\tag 2.30
$$ 
with strongly homogeneous symbols $s_{J,j,j'}$ of degree $-3-J+j+j'$.

We have shown:

\proclaim{Lemma 2.4} The singular Green operator $G_\mu $ in $R_\mu
=Q_{\mu ,+}+G_\mu $ has the
form {\rm (2.28)},
with strongly polyhomogeneous entries.
In local coordinates, the symbols of the Poisson operators $K_{0,\mu
}$ and $K_{1,\mu }$ have 
homogeneous terms that are rational functions of $\xi _n$ with the
pole $i\kappa ^+(x',\xi ',\mu )$, and the symbol of the trace
operator $T_{0,\mu }$ 
has homogeneous terms that are rational functions of $\xi _n$ with the
pole $-i\kappa ^-(x',\xi ',\mu )$, as described in detail in {\rm (2.29)ff}.
The symbol structure of $G_\mu $ is described in \rm{(2.30)ff}.
\endproclaim

Moreover, $G^{(k)}_\mu $ has the same kind of symbol
structure as $G_\mu $:

\proclaim{Proposition 2.5} $G^{(k)}_\mu $ in {\rm (2.3)} is a strongly
polyhomogeneous 
singular Green operator of order $-2k$ (degree $-1-2k$), and the homogeneous terms
$g^{(k)}_{-1-2k-J}(x',\xi ',\xi _n,\eta _n,\mu )$ 
in its symbol $g^{(k)}\sim \sum_{J\in\Bbb
N}g^{(k)}_{-1-2k-J}$ 
are of the form $$
g^{(k)}_{-1-2k-J}(x',\xi ',\xi _n,\eta _n,\mu )=\sum_{ \Sb j\ge 1,\,j'\ge 1\\
j+j'\le 2J+k+1\endSb}\tfrac
{s_{k,J,j,j'}(x',\xi ',\mu )}{(\kappa ^+(x',\xi ',\mu )+ i\xi
_n)^{j}(\kappa ^-(x',\xi ',\mu )- i\eta 
_n)^{j'}},\tag2.31
$$
with $s_{k,J,j,j'}$ strongly homogeneous of degree $-1-2k-J+j+j'$. 
\endproclaim 

\demo{Proof} For the term $(Q_{\mu ,+})^k-(Q_{\mu }^k)_+$ we have:
$$\aligned
(Q_{\mu ,+})^k-(Q_{\mu }^k)_+
&=(Q_{\mu ,+})^k-(Q_{\mu }^2)_+(Q_{\mu
,+})^{k-2}+ (Q_{\mu }^2)_+(Q_{\mu
,+})^{k-2}-(Q_{\mu }^3)_+(Q_{\mu
,+})^{k-3}\\
&\quad +(Q_{\mu }^3)_+(Q_{\mu
,+})^{k-3}+\dots +(Q_{\mu }^{k-1})_+Q_{\mu
,+}-  
(Q_{\mu }^k)_+\\
&=-G^+(Q_\mu )G^-(Q_\mu )(Q_{\mu ,+})^{k-2}-G^+(Q^2_\mu )G^-(Q_\mu )(Q_{\mu
,+})^{k-3}\\
&\quad-\dots-
G^+(Q^{k-1}_\mu )G^-(Q_\mu ).
\endaligned\tag2.32$$
By Proposition 2.3, the terms in the symbol of $G^+(Q^{k'})G^-(Q_\mu )$ are
compositions of the form$$
\tfrac{r_1}{(\kappa ^++i\xi _n)^{j_1}(\kappa ^+-i\eta _n)^{j_2}}\circ \tfrac{r_2}{(\kappa ^-+i\xi _n)^{j_3}(\kappa ^--i\eta _n)^{j_4}}.
$$
The composition rule in the normal variable is recalled in (3.7)
below; the composition in tangential variables follows standard $\psi
$do rules. As a result, we get sums of expressions $$
\tfrac{r_3}{(\kappa ^++i\xi _n)^{j_5}(\kappa ^--i\eta _n)^{j_6}},\tag2.33
$$ since 
$\tfrac1{(\kappa ^+-i\eta _n)^{j_2}}\circ_n\tfrac1{(\kappa
^-+i\xi _n)^{j_3}}=\frac1{2\pi }\int \tfrac1{(\kappa ^+-i\xi _n)^{j_2}}\tfrac1{(\kappa
^-+i\xi _n)^{j_3}}\, d\xi _n$ gives a strongly polyhomogeneous $\psi
$do symbol on $\Bbb R^{n-1}$.
Next, when an s.g.o\. with symbol of the form (2.33) is composed with
$Q_{\mu ,+}$, we must calculate$$
h^-_{\eta _n}\Bigl[\tfrac{1}{(\kappa ^++i\xi _n)^{j_5}}\tfrac1{(\kappa
^--i\eta _n)^{j_6}}\Bigl( \tfrac{1}{(\kappa ^++i\eta _n)^{j_7}}+\tfrac1{
(\kappa ^--i\eta _n)^{j_7}}\Bigr)\Bigr]\tag2.34
$$
(as recalled in (3.9) below), which is seen to give symbols of the
form (2.33) by decomposition of 
the function of $\eta _n$ in
simple fractions. So indeed, the contribution from (2.32) has
homogeneous terms of the form (2.33). 

For the composition $G_\mu Q_{\mu ,+}$ we have the same situation as
in (2.34), resulting in another s.g.o\. of the asserted type. For
compositions $Q_{\mu ,+}G_\mu $ we have instead to calculate terms of
the form$$
h^+_{\xi _n}\Bigl[\Bigl( \tfrac{1}{(\kappa ^++i\xi  _n)^{j_7}}+\tfrac1{
(\kappa ^--i\xi  _n)^{j_7}}\Bigr)\tfrac{1}{(\kappa ^++i\xi _n)^{j_5}}\tfrac1{(\kappa
^--i\eta _n)^{j_6}}\Bigr],
$$
which likewise lead to expressions as in (2.33). Further compositions
behave similarly. Thus also the
contributions from $\operatorname{pol}(G_\mu ,Q_{\mu ,+})$ are of the
asserted form.
\qed
\enddemo 

\subhead 3. Trace calculations for terms in $G^+(P)G^-(Q_\mu ^k)$
and $GG_\mu ^{(k)}$  
\endsubhead

In the next three sections we analyze the traces of the operators
$G^+(P)G^-(Q^k_\mu )$,
$P_+G^{(k)}_\mu$,
$GQ^k_{\mu ,+}$ and 
$GG^{(k)}_\mu$
in the localized situation. We consider in detail the case where the
fiber dimension $n'$ of $E$ equals 1; it serves as a model, and the
results in general are found by adding finitely many terms of this
kind.
For the right hand side
factor in each composition we insert 
the expansion of the symbol in homogeneous terms. The 
remainders $
G^-(Q_\mu ^k)-\operatorname{OPG}(\sum_{J<J_0}g^{\circ k,-}_{-1-2k-J})$,
$G^{(k)}_\mu -\operatorname{OPG}(\sum_{J<J_0}g^{( k)}_{-1-2k-J})$,
$Q_{\mu ,+}^k-\operatorname{OP}(\sum_{J<J_0}q^{\circ
k}_{-2k-J})_+$,
are bounded from $L_2(\rnp)$ to $H^{2k+J_0,\mu }(\crnp)$ uniformly for
$\mu $ on rays (cf\. [G96, Theorems 2.5.6, 2.5.11]), so the
compositions of these with $G^+(P)$, $P_+$ or $G$ have norm $O(\mu
^{-J_0})$ as operators from  $L_2(\rnp)$ to the space of functions
in $H^{2k-\nu }(\crnp)$ with support in a fixed compact set. Since
$2k-\nu >n$, they are trace class with trace norms that are $O(\mu
^{-J_0})$:$$\multline
\|G^+(P)[G^-(Q_\mu ^k)-\operatorname{OPG}(\sum_{J<J_0}g^{\circ
k,-}_{-1-2k-J})]\|_{\Tr},\quad 
\|P_+[G^{(k)}_\mu -\operatorname{OPG}(\sum_{J<J_0}g^{(
k)}_{-1-2k-J})]\|_{\Tr},\\
\|G[G^{(k)}_\mu -\operatorname{OPG}(\sum_{J<J_0}g^{(
k)}_{-1-2k-J})]\|_{\Tr},\quad
\|G[
Q_{\mu ,+}^k-\operatorname{OP}(\sum_{J<J_0}q^{\circ
k}_{-2k-J})_+]\|_{\Tr}\\
\text{ are all }O(\mu ^{-J_0}).\endmultline\tag3.1
$$
Here we can take $J_0$ arbitrarily large, so it remains to
analyze the traces of
\linebreak$G^+(P)\operatorname{OPG}(g^{\circ
k,-}_{-1-2k-J})$,
$P_+\operatorname{OPG}(g^{(
k)}_{-1-2k-J})$,
$G\operatorname{OPG}(g^{(
k)}_{-1-2k-J})$ and
$G\operatorname{OP}(q^{\circ
k}_{-2k-J})_+$,
for each $J$.

In the following, we operate with symbols in $(x',y')$-form and in
$y'$-form as 
well as $x'$-form (cf\. [G96], p\. 141, Sect\. 2.4]); the advantage is that 
$$g_1(x',\xi ,\eta _n )\circ g_2(y',\xi ,\eta _n)=g_1(x',\xi ,\eta _n
)\circ_n g_2(y',\xi ,\eta _n),
\tag3.2
$$ with no further terms; this may be
reduced to $x'$-form by [G96, Th\. 2.4.6].

We shall use the fact that a singular Green operator
$G=\operatorname{OPG}(g(x',y',\xi ',\xi _n,\eta 
_n))$ on $\rnp$  
of order $<1-n$ and class 0 with compact $(x',y')$-support is
trace class (in view of the rapid decrease of the symbol-kernel
$\tilde g$ (cf\. (2.18)) for
$x_n,y_n\to\infty $), and that its trace can be
calculated by a diagonal integral of its  
kernel $
\Cal K_{G}(x,y)$:$$
\aligned
 \Cal K_{G}(x,y)&=(2\pi )^{-n-1}\int _{\Bbb
R^{n+1}}e^{i(x'-y')\cdot\xi 
'+ix_n\xi _n-iy_n\eta _n}g(x',y',\xi ',\xi _n ,\eta _n)\,d\xi'd\xi _n d\eta _n\\
&=(2\pi )^{1-n}\int _{\Bbb
R^{n-1}}e^{i(x'-y')\cdot\xi 
'}\tilde g(x',y',x_n,y_n,\xi ')\,d\xi ',\\
\Tr G &=\int_{\rnp}\Cal
K_{ G}(x,x)\,dx=(2\pi )^{1-n}\int_{\rnp\times\Bbb R^{n-1}}\tilde
g(x',x',x_n,x_n,\xi ')\,d\xi 'dx.
\endaligned\tag 3.3
$$
(This and the following holds also when $G$ is in $x'$-form or
$y'$-form, vanishing for large $x'$ resp\. $y'$.)
One defines the normal trace $\tr_ng$ by:
$$ \tr_n g(x',y',\xi
')=\int_0^\infty \tilde g(x',y',x_n,x_n,\xi ')\,dx_n=\tfrac1{2\pi }\int g(x',y',\xi ',\xi
_n,\xi _n)\,d\xi _n,
\tag 3.4
$$it is a $\psi $do symbol on $\Bbb R^{n-1}$. ($\tr_n g(x',\xi ')$ is called
$\ttilde g(x',\xi ')$ in [G92], [G96] and \linebreak$\tr g(x',\xi ',D_n)$ in
[FGLS96].) Denoting by $\tr_nG$ the
$\psi $do it defines, 
we obtain for its trace:
$$
\Tr'(\tr_n G)
=(2\pi )^{n-1}\int_{\Bbb R^{2n-2}}\int_0^\infty 
\tilde
g(x',x',x_n,x_n,\xi ')\,dx_n d\xi 'dx'=\Tr G.
\tag 3.5$$
Here we write $\Tr'$ to indicate the trace of an
operator defined over $\Bbb R^{n-1}$.

When $g(x',y',\xi ,\eta _n,\mu
)$ is strongly polyhomogeneous and has order $\nu <1-n$
and compact $(x',y')$-support (or is in $x'$-form or $y'$-form with
compact $x'$-support resp\. $y'$-support), then the
operator it defines is trace class with an expansion
$$
\Tr \operatorname{OPG}(g)\sim \sum_{r\in\Bbb N}c_r\mu ^{n-1+\nu -r}\text{
for }\mu \to\infty \text{ in }\Gamma ,\tag 3.6
$$
without logarithmic terms. This holds by the proof of [G96, Th\.
3.3.10ff.] and its 1986 predecessor, also recalled in [G92, App.]; the
reason is simply that $\tr_n g$ (or $\ttilde g$) is a $\psi $do symbol
of order $\nu $ and regularity $+\infty $ so that the homogeneous terms
in its symbol give the terms in (3.6).

\example
{Example 3.1} To explain the basic point of our strategy, we consider
the simple case where $P=0$ and $G$ and $P_1$ are replaced by operators with
$x$-independent symbols, $G=\operatorname{OPG}(g(\xi ,\eta _n))$ and
$P_1=1-\Delta =\operatorname{OP}(p_1(\xi ))$, with
$$
g(\xi ,\eta _n)=\frac{c(\xi ')}{(\sigma (\xi ')+i\xi
_n)(\sigma (\xi ')-i\eta  _n)},\quad p_1(\xi )=|\xi |^2+1,\quad
\sigma (\xi ')=[\xi '];$$
here $c(\xi ')$ is a $\psi $do symbol on $\Bbb R^{n-1}$ of order $\nu
+1$, whereby $G$ is of order $\nu $; 
$[\xi ']$ denotes a positive $C^\infty $ function equal to
$|\xi '|$ for $|\xi '
|\ge \frac12$. We assume $\nu -2<-n$. 

The resolvent $R_\mu =
(P_{1,D}+\mu^2)^{-1}$  is of the form $Q_{\mu,+}+G_{1,\mu}$, where  
$G_{1,\mu} = -K_\mu\gamma_0 Q_{\mu,+}$ with $K_\mu$ denoting the
Poisson operator $K_\mu \:v\mapsto u$ 
solving the Dirichlet problem \linebreak
$(P_1+\mu^2)u = 0$, $\gamma_0 u=v$. 
It is easily checked that $K_\mu$ has the symbol-kernel 
$e^{-x_n\kappa}$ and symbol
$(\kappa+i\xi_n)^{-1}$, with $\kappa = \langle (\xi',\mu)\rangle$.
Thus (cf\. [G96, Theorem 2.6.1] for composition rules), 
 $R_\mu = \operatorname{OP}(q)_+
+\operatorname{OPG}(g_1)$ with  symbols
$$\aligned
q(\xi ,\mu
)&=\frac1{|\xi |^2+1+\mu ^2}=\frac1{(\kappa +i\xi _n)(\kappa -i\xi
_n)}=
\frac1{2\kappa (\kappa +i\xi _n)}+\frac1{2\kappa (\kappa -i\xi
_n)},\\
g_1(\xi ,\eta _n,\mu
)&=-\frac1{\kappa +i\xi _n}h^-q(\xi ',\eta _n,\mu )=\frac{-1}{2\kappa (\kappa +i\xi _n)(\kappa -i\eta 
_n)} .\endaligned
$$ 

Consider $GG_\mu $. Composition with respect to the normal
variable (in fact the full composition since the symbols are independent of
$x'$) and application of $\tr_n$ have the form
$$\gathered
g(\xi ',\xi_n, \eta _n)\circ_n g_1(\xi ',\xi _n,\eta _n,\mu
)=\tfrac1{2\pi }\int 
g(\xi',\xi _n ,\zeta _n)g_1(\xi ',\zeta _n,\eta _n,\mu )\,d\zeta
_n,\\
\tr_n (g\circ_n g_1)=\tfrac1{(2\pi )^2}\int g(\xi',\xi _n ,\zeta
_n)g_1(\xi ',\zeta _n,\xi  _n,\mu )\,d\zeta _nd\xi _n,
\endgathered
\tag3.7 
$$ 
which in the considered case gives$$
\tr_n (g\circ_n g_1)=\tfrac1{(2\pi )^2}\int \frac{c(\xi ')}{(\sigma+i\xi
_n)(\sigma-i\zeta  _n)}\frac{-1}{2\kappa (\kappa +i\zeta
_n)(\kappa -i\xi  
_n)}\,d\zeta _nd\xi _n=\frac{-c}{2\kappa (\sigma+\kappa )^2},
$$
by residue calculus. This is a product of $\psi $do symbols on $\Bbb
R^{n-1}$ that are weakly
polyhomogeneous according to [GS95],  with  $c\in S^{\nu +1,0}$,
$\kappa ^{-1}$ and 
$(\sigma+\kappa )^{-1}\in S^{-1,0}\cap S^{0,-1}$, so $\tr_n (g\circ_n
g_1)\in S^{\nu -2,0}\cap S^{\nu +1,-3}$. 
When $\psi (x')\in C_0^\infty (\Bbb R^{n-1})$, the operator $
\psi \operatorname{OP}'(\tr_n 
g\circ g_1)=\psi \tr_n(GG_\mu ) $ is trace
class since $\nu -2<-n$ (in fact $\nu -2<1-n$
would suffice). By (3.5) and
[GS95, Th\. 2.1],$$
\Tr'(\psi \tr_n GG_\mu )=\Tr(\psi GG_\mu ) \sim
\sum_{l=0}^\infty a_{l}\mu ^{n-1+\nu
-2-l}+\sum_{l=0}^\infty (a'_l\log\mu +a''_l)\mu ^{-3-l};\tag3.8
$$
note that there is {\it no log-term at the
power }$-2$. 

Now consider $GQ_{\mu ,+}$.
The composition rules are here 
$$\gathered
g(\xi ',\xi_n, \eta _n)\circ_n q(\xi ',\xi _n,\mu
)_+=h^-_{\eta _n}[g(\xi',\xi _n ,\eta _n)q(\xi ',\eta _n,\mu )],\\
\tr_n(g\circ_n q_+)
=\tfrac1{2\pi }\int 
h^-_{\eta _n}[g(\xi',\xi _n ,\eta _n)q(\xi ',\eta _n,\mu )]|_{\xi
_n=\eta _n}\,d\eta _n\\
=\tfrac1{2\pi }\int 
g(\xi',\eta _n ,\eta _n)q(\xi ',\eta  _n,\mu )\,d\eta _n,
\endgathered\tag3.9 
$$ 
(since $h^+_{\eta _n}[g(\xi ',\xi _n,\eta _n)q(\xi ',\eta _n,\mu
)]|_{\xi _n=\eta _n}$
is meromorphic in $\eta _n$ with no poles in $\Bbb C_-$ and is
$O(\eta _n^{-2})$
for $\eta _n\to\infty $ in $\Bbb C$).
In the considered case this gives:$$
\tr_n (g\circ_n q_+)=\tfrac1{(2\pi )^2}\int \frac{c(\xi ')}{(\sigma+i\xi
_n)(\sigma-i\xi _n)}\Bigl[\frac{1}{2\kappa  (\kappa +i\xi
_n)}+\frac1{2\kappa (\kappa -i\xi 
_n)}\Bigr]\,d\xi _n=\frac{c}{2\kappa (\sigma+\kappa )\sigma},
$$
by residue calculus. Then since  $\sigma^{-1}\in S^{-1,0}$, we find
that
$
\tr_n (g\circ_n q_{+})\in S^{\nu -2,0}\cap S^{\nu ,-2}
$.
By [GS95, Th\. 2.1],
$$
\Tr'(\psi \tr_n GQ_{\mu,+}) =\Tr(\psi GQ_{\mu,+}) \sim
\sum_{l=0}^\infty b_{l}\mu ^{n-1+\nu
-2-l}+\sum_{l=0}^\infty (b'_l\log\mu +b''_l)\mu ^{-2-l}.\tag3.10
$$
In this case there {\it is} a log-term at the power $-2$. 
The proof of [GS95, Th\. 2.1] allows an
analysis of 
the coefficient $b'_0$ in (3.10).
It shows that $b'_0= (2\pi)^{1-n}\int_{S^*X'}
s(x',\xi')_{1-n} d\sigma'$,
where $s(x',\xi')$ is the coefficient of $\mu^{-2}$ in the expansion
of $\psi \tr_n (g\circ_n q_{+})$ into decreasing powers of $\mu$, and
the subscript denotes the component of homogeneity $1-n$ with respect to 
$\xi'$.
In our case, we have  $
\psi \tr_n (g\circ_n q_{+})=\psi c(2\kappa (\kappa +\sigma )\sigma)^{-1}\sim
\psi c(2\sigma)^{-1}\mu ^{-2}+\ldots$ (lower order in $\mu$). Since $\tr_n
(\psi g)=\psi c(2\sigma)^{-1}$, we obtain the statement for
$\psi G$ on the
noncommutative residue in Theorem 1.1.
\endexample

In the general case we shall use expansions of the
$\mu $-independent symbols in Laguerre functions (with the pole $
i\sigma )$, which together with the structure of the
symbols of $Q_\mu $ and $G_\mu $ as rational functions with
well-behaved poles allow termwise
residue calculus like in the example. The general outcome is that the
compositions with $\mu $-dependent s.g.o.s do not contribute to the relevant
logarithmic term, whereas the compositions with $Q_{\mu ,+} ^k$ do so. It
is important to show that the
results we find are preserved after summation of the Laguerre expansions.
The reader is kindly advised to consult
the Appendix, where we have collected some facts on the Laguerre
expansions used here.

In the present section we study the terms in $G^+(P)G^-(Q_\mu ^k)$
and $GG_\mu ^{(k)}$.
These two compositions 
are treated in similar ways; we give the details for
the latter which contains the most general types of singular Green operators.

In the local coordinate system, $G=\operatorname{OPG}(g(x',\xi ',\xi
_n,\eta _n))$, where we can expand
the symbol $g(x',\xi ',\xi _n,\eta _n)$ in a convergent series 
in terms of Laguerre functions of
$\xi _n$ and $\eta _n$:
$$\aligned
g(x',\xi ',\xi _n,\eta _n)&=\sum_{l,m\in\Bbb N}2\sigma c_{lm}(x',\xi
')\hat\varphi '_l(\xi _n,\sigma )
\bar{\hat\varphi }'_m(\eta _n,\sigma ),\quad \sigma =[{\xi '}],\\
\endaligned\tag3.11
$$
with a rapidly decreasing double sequence of $\psi $do symbols
$c_{lm}(x',\xi ')$ of order $\nu $ with compact $x'$-support.
(Note that we use the non-normed Laguerre functions $\hat\varphi
'_l$, cf\. (A.4).) For $g^{(
k)}_{-1-2k-J}$ we have the structure
described in Proposition 2.5.

We shall use formulas (3.4)--(3.5) to calculate the
trace of $G\,\operatorname{OPG}(g^{(
k)}_{-1-2k-J})$. For convenient composition rules,
we assume that the right hand factor is given in $y'$-form, with symbol
$g^{(k)}_{-1-2k-J}(y',\xi ,\eta _n,\mu )$ (recall (3.2)). 
The general rule for passage from $x'$-form to $y'$-form (cf\.
[G96, (2.4.30 iii)]) shows 
that this
symbol has the same structure as the $x'$-form described in
Proposition 2.5. 
Then$$\aligned
\Tr
G\,\operatorname{OPG}(g^{(k)}_{-1-2k-J})&=\Tr\operatorname{OPG}(g(x',\xi ,
\eta _n)\circ_n 
g^{(k)}_{-1-2k-J}(y',\xi ,\eta _n,\mu )),\text{ where}\\
g\circ_n g^{(k)}_{-1-2k-J}&=\sum_{l,m\in\Bbb
N}g_{l,m,J}(x',y',\xi ,\eta _n,\mu )\text{ with}\\
g_{l,m,J}&=
2\sigma c_{lm}(x',\xi
')\hat\varphi '_l(\xi _n,\sigma )
\bar{\hat\varphi }'_m(\eta _n,\sigma )\circ_n 
g^{(k)}_{-1-2k-J}(y',\xi ,\eta _n,\mu ).
\endaligned\tag3.12
$$
Termwise, we have by (3.5),$$\aligned
\Tr\operatorname{OPG}(g_{l,m,J}(x',y',\xi ,\eta _n,\mu
))&=\Tr'\operatorname{OP}'( s_{l,m,J}(x',y',\xi ',\mu )),\text{
where}\\
 s_{l,m,J}(x',y',\xi ',\mu )&= \tr_ng_{l,m,J}(x',y',\xi ,\eta _n,\mu
).
\endaligned\tag 3.13
$$
The crucial
step is now to calculate the normal trace of the 
s.g.o\. arising from composing $\hat\varphi '_l(\xi _n,\sigma )
\bar{\hat\varphi }'_m(\eta _n,\sigma )$ with the fractions in (2.31).
In view of (3.7), we must calculate
$$\aligned
\tr_n &\bigl( \hat\varphi '_l(\xi _n,\sigma )\bar{\hat\varphi
}'_m(\eta _n,\sigma )\circ_n\tfrac1{(\kappa ^++i\xi _n)^j(\kappa ^--i\eta
_n)^{j'}}\bigr)\\
&=\tr_n \bigl( \hat\varphi '_l(\xi  _n,\sigma )\bigl[\tfrac1{2\pi }\tsize\int \bar{\hat\varphi
}'_m(\zeta  _n,\sigma )\tfrac1{(\kappa ^++i\zeta
_n)^j}\,d\zeta _n \bigr]\tfrac 1{(\kappa ^--i\eta  
_n)^{j'}}\bigr)\\
&=\tfrac1{2\pi }\tsize\int \bar{\hat\varphi
}'_m(\xi _n,\sigma ) \frac1{(\kappa ^++i\xi _n)^j}\,d\xi _n\cdot
\tfrac1{2\pi }\tsize\int \hat\varphi '_l(\xi  _n,\sigma )\tfrac
1{(\kappa ^--i\xi  
_n)^{j'}} \,d\xi _n; 
\endaligned\tag3.14
$$
a product of two $\psi $do symbols on $\Bbb R^{n-1}$.

The resulting formulas are worked out in the following lemma. 

\proclaim{Lemma 3.2} One has for all $m\ge 0$ and $j\ge 1$:
$$\aligned
\tfrac1{2\pi }\int \bar{\hat\varphi }'_{m}(\xi _n,\sigma )\frac{1}
{(\kappa ^+(x',\xi ',\mu )+i\xi
_n)^{j}}\,d\xi _n&=\sum_{m'\ge 0,\, |m'-m|<j}(\kappa^+) ^{1-j} a_{jm'}\frac{(\sigma -\kappa ^+)^{m'}}{(\sigma
+\kappa ^+ 
)^{m'+1}},\\
\tfrac1{2\pi }\int {\hat\varphi }'_{m}(\xi _n,\sigma )\frac{1}
{(\kappa ^-(x',\xi ',\mu )-i\xi
_n)^{{j}}}\,d\xi _n
&=\sum_{m'\ge 0,\, |m'-m|<j}(\kappa ^-) ^{1-{j}} \bar a_{{j}m'}\frac{(\sigma -\kappa ^-)^{m'}}{(\sigma +\kappa ^-
)^{m'+1}},
\endaligned\tag3.15$$
with universal constants $a_{jm'}$ that are $O(m^j)$ for fixed $j$. The
resulting expressions are
weakly polyhomogeneous $\psi $do symbols belonging to $S^{-j,0}(\Bbb R^{n-1}\times\Bbb R^{n-1},\Gamma )\cap
S^{0,-j}(\Bbb R^{n-1}\times\Bbb R^{n-1},\Gamma )$.
\endproclaim 

\demo{Proof} 
If $j=1$, then since the integrand is $O(\xi _n^{-2})$ for $|\xi
_n|\to\infty $,
$$\multline
\tfrac1{2\pi }\int \bar{\hat\varphi }'_{m}(\xi _n,\sigma )\frac{1}
{\kappa ^+(x',\xi ',\mu )+i\xi
_n}\,d\xi _n
=\tfrac1{2\pi i
}\int_{\Cal L_+}\frac{(\sigma +i\xi _n)^{m}}{(\sigma -i\xi
_n)^{m+1}}\frac1{\xi _n-i\kappa ^+}\,d\xi _n\\
=\operatorname{Res}_{\{\xi
_n=i\kappa ^+\}}\frac{(\sigma +i\xi _n)^{m}}{(\sigma -i\xi
_n)^{m+1}}\frac1{\xi _n-i\kappa ^+}
=\frac{(\sigma +i\xi _n)^{m}}{(\sigma -i\xi
_n)^{m+1}}\Big|_{\xi _n=i\kappa ^+}=\frac{(\sigma -\kappa
^+)^{m}}{(\sigma +\kappa ^+ 
)^{m+1}},\endmultline\tag3.16
$$where $\Cal L_+$ is a positively oriented curve around $i\kappa ^+$ in
$\Bbb C_+$. By the rules in [GS95] (more specifically, Th\.
1.16 and 1.23 there), $\frac{(\sigma -\kappa 
^+)^{m}}{(\sigma +\kappa ^+
)^{m+1}}$ is weakly polyhomogeneous (wphg) belonging to $S^{-1,0}
(\Bbb R^{n-1}\times\Bbb R^{n-1},\Gamma )\cap
S^{0,-1}(\Bbb R^{n-1}\times\Bbb R^{n-1},\Gamma )$ (recall (2.7)).
For a larger value of $j$ we use (A.7)
and the third line in (A.5) (conjugated), noting that $$
\xi _n^{j-1}\partial
_{\xi _n}^{j-1}f(\xi _n)=\sum_{0\le j'\le j-1}c_{jj'}(\partial _{\xi
_n}\xi _n)^{j'}f(\xi _n),\tag3.17
$$
with universal constants $c_{jj'}$.
This gives the following calculation:
$$\aligned
\tfrac1{2\pi
}\int_{\Cal L_+}\bar{\hat\varphi }'_{m}&(\xi _n,\sigma )\frac1{(\kappa
^++i\xi _n)^j}\,d\xi _n=\tfrac1{2\pi
}\int_{\Cal L_+}\bar{\hat\varphi }'_{m}(\xi _n,\sigma
)\tfrac1{(j-1)!}(i\partial _{\xi _n})^{j-1}\frac1{\kappa ^++i\xi 
_n}\,d\xi _n\\&=\tfrac1{2\pi
}\int_{\Cal L_+}\tfrac1{(j-1)!}[(-i\partial _{\xi
_n})^{j-1}\bar{\hat\varphi }'_{m}(\xi _n,\sigma )]\frac1{\kappa ^++i\xi
_n}\,d\xi _n\\ 
&=\tfrac{(-i)^{j-1}}{(j-1)!}
\partial _{\xi _n}^{j-1}\bar{\hat\varphi }'_{m}(\xi _n,\sigma
)\Big|_{\xi _n=i\kappa ^+}=\tfrac{(-\kappa ^+)^{1-j}}{(j-1)!}
\xi _n^{j-1}\partial _{\xi _n}^{j-1}\bar{\hat\varphi }'_{m}(\xi
_n,\sigma )\Big|_{\xi _n=i\kappa ^+}\\
\endaligned$$
$$\aligned
&=\tfrac{(-\kappa ^+)^{1-j}}{(j-1)!}
 \sum_{0\le j'\le j-1}c_{jj'}(\partial _{\xi _n}\xi _n)^{j'}\bar{\hat\varphi }'_{m}(\xi
_n,\sigma )\Big|_{\xi _n=i\kappa ^+}\\
&= \sum_{m'\ge 0,\, |m'-m|<j}{(\kappa ^+)^{1-j}}a_{jm'}\bar{\hat\varphi }'_{m'}(\xi _n,\sigma )\Big|_{\xi _n=i\kappa ^+}\\
&= \sum_{m'\ge 0,\, |m'-m|<j}(\kappa^+) ^{1-j}a_{jm'}\frac{(\sigma -\kappa ^+)^{m'}}{(\sigma +\kappa ^+
)^{m'+1}},
\endaligned\tag3.18 $$
with some universal constants $a_{jm'}$ that are $O(m^j)$. (Observe
that in the third line of (A.5), the fact that $\hat\varphi '_{k-1}$
in the right hand side has coefficient
$-\frac12 k$ assures that no terms with
negative index are introduced when we apply the formula repeatedly to a
$\hat\varphi '_m$ with $m\ge 0$.) In view of the observations
made further above on $\frac{(\sigma -\kappa ^+)^{m}}{(\sigma +\kappa ^+
)^{m+1}}$, and the composition rules of [GS95], this is a
wphg $\psi $do symbol in $S^{-j,0}\cap
S^{0,-j}$.
This proves the statements for
the first line in (3.15).

The second line in (3.15) follows by conjugation, replacing $\bar\kappa ^+$ by
$\kappa ^-$.
\qed
\enddemo 

\example{Remark 3.3} One has moreover that the
$\partial _\mu ^r$-derivatives belong to $S^{-j-r,0}\cap S^{0,-j-r}$ for
all $r$ (such symbols are called special
parameter-dependent of 
degree $-j$ in [G99]). Similar statements hold for symbols
constructed in the following. 
\endexample

We then find:

\proclaim{Proposition 3.4} 
The pseudodifferential symbols $
 s_{l,m,J}(x',y',\xi ',\mu )$
defined in {\rm (3.13)}
are in $S^{\nu -2k-J,0}(\Bbb R^{n-1}\times\Bbb R^{n-1},\Gamma )\cap
S^{\nu +1,-1 -2k-J}(\Bbb R^{n-1}\times\Bbb R^{n-1},\Gamma )$.

Moreover,
the sum over $l,m$ converges in all symbol norms, so that also the
symbol$$
\multline
s_{(J)}(x',y',\xi ',\mu )=\sum_{l,m\in\Bbb N}s_{l,m,J}(x',y',\xi
',\mu )\\
=\tr_n(g(x',\xi ',\xi _n,\eta _n)\circ_n g^{(k)}_{-1-2k-J}(y',\xi ',x_n,\eta _n,\mu ))\endmultline\tag3.19
$$
is in 
$S^{\nu -2k-J,0}\cap
S^{\nu +1,-1-2k-J}$.
 
Finally, there is an expansion for $\mu \to \infty $ in $\Gamma $:
$$\multline
\Tr\bigl(G\,\operatorname{OPG}(
 g^{(k)}_{-1-2k-J}(y',\xi ,\eta _n,\mu
))\bigr)
= \Tr'\operatorname{OP}'\bigl(s_{(J)}(x',y',\xi ',\mu )\bigr)\\
\sim \sum_{l=0}^\infty a^{1}_{J,l}\mu ^{n-1+\nu
-2k-J-l}+\sum_{l=0}^\infty (a^{2}_{J,l}\log\mu +a^{3}_{J,l})\mu ^{-1 -2k-J-l}.
\endmultline \tag3.20
$$

\endproclaim

\demo{Proof} 
The contribution from each term in the finite expansion (2.31) is
treated by Lemma 3.2. For 
$s_{k,J,j,j'}(x',\xi ',\mu )$ and $2\sigma c_{lm}$ replaced by 1, this
gives a
symbol in $S^{-j-j',0}\cap
S^{0,-j-j'}$. The
factors are taken into account as follows: Since 
$s_{k,J,j,j'}(x',\xi ',\mu )$ is strongly homogeneous of degree
$-1-2k-J+j+j'$, it is 
in $S^{-1-2k-J+j+j',0}\cap
S^{0,-1-2k-J+j+j'}$
so multiplication by it gives a
symbol in $S^{-1-2k-J,0}\cap
S^{0,-1-2k-J}$. Next,
 $2\sigma c_{lm}\in S^{\nu +1,0}$, so multiplication by it gives a
symbol in $S^{\nu -2k-J,0}\cap
S^{\nu +1,-1-2k-J}$.
This shows the statement on the $ s_{m,l,J}$.

For the summation in $l,m$, we use that the symbols $c_{lm}(x',\xi
')$ are rapidly decreasing in $l$ 
and $m$ in all 
symbol norms. Observe moreover, that the new coefficients
introduced in the formulas in 
(3.15) are polynomially bounded in $l$ and $m$, for each $j$ and $j'$,
and that the same will hold for derivatives in view of the formulas 
(with $\kappa =\kappa ^\pm$),$$\aligned
\partial _{\xi _j}(\tfrac {\sigma -\kappa }{\sigma +\kappa })^{m}&=
m(\tfrac {\sigma -\kappa }{\sigma +\kappa })^{m-1}\tfrac{\kappa \partial
_{\xi _j}\sigma  -\sigma \partial _{\xi _j}\kappa }{(\sigma
+\kappa )^2},\\
\partial _{z}(\tfrac {\sigma -\kappa }{\sigma +\kappa })^{m}&=
m(\tfrac {\sigma -\kappa }{\sigma +\kappa })^{m-1}\tfrac{ -\sigma
\partial _{z}\kappa }{(\sigma +\kappa )^2};
\endaligned\tag3.21
$$ that can be continued to a study of higher derivatives
again giving expressions with polynomial coefficients in
$m$. Using these facts together, one finds that
summation in $l,m$ converges in each symbol norm. This shows the
second assertion.

Since the trace norm is dominated by a finite set of symbol norms,
the summation now also converges in trace norm, so by (3.13),$$
\Tr'\operatorname{OP}'(s_{(J)})=\sum_{l,m}\Tr'\operatorname{OP}'(
s_{l,m,J})=\sum_{l,m}\Tr\operatorname{OPG}(g_{l,m,J})=\Tr\operatorname{OPG}(g\circ_n
g^{(k)}_{-1-2k-J});
$$
for the last equality sign we used that the sum of s.g.o.s likewise
converges in trace norm. Now we can apply [GS95, Th\. 2.1] to
$\operatorname{OP}'(s_{(J)})$ in $n-1$ dimensions, which shows that the trace 
has an asymptotic expansion (3.20).\qed
\enddemo

Note that the expansions (3.20) are of the type (1.10) (I) with $-\lambda
=\mu ^2$, and that they {\it do not contribute to the coefficient $c'_0$
there,} not even for $J=0$,
since the first logarithmic term in (3.20) has $\mu ^2$ to the power
$-k-\frac12-\frac J2$.

A similar result holds for the operators
$G^+(P)\,\operatorname{OPG}(g^{\circ k,-}_{-1-2k-J})$, except that
the root in 
$\Bbb C^+$ is now $i\kappa ^-$; we shall not write a separate
statement on this.

\example{Remark 3.5}
One can fairly easily allow $P_1$ to be of higher order (or non-scalar)
with constant multiplicity of the roots (in $\xi _n$) of the
(determinant of the) principal symbol at $x_n=0$, since the symbol of
$G^{(k)}_\mu $ then 
has a decomposition in rational functions in terms of these poles
depending continuously on $(\xi ',\mu )$; this just requires a
little more residue calculus (with several poles instead of one in
$\Bbb C_\pm$).
 For the general case where multiplicities may vary with $(\xi ',\mu )$, one
can probably instead obtain the desired results by use of contour
integrals around all the poles in $\Bbb C_{\pm}$, respectively. 
We leave out these generalizations here, since the method is already
complicated to explain in the case of the two poles $\kappa ^+$ and
$\kappa ^-$.

\endexample

\subhead 4. Trace calculations for terms in $P_+G_\mu ^{(k)}$ \endsubhead

We now turn to the operator $P_+G_\mu ^{(k)}$ in (2.5), in the
localized version. Here $P_+G_\mu ^{(k)}=P_+G_\mu ^{(k)}\theta $ for
some $\theta \in C_{(0)}^\infty (\crnp)$ ($C^\infty $ functions with
compact support in $\crnp$). Some manipulations with the effect of
$P$ are needed before we decompose the right hand factor in
homogeneous terms.
We can assume that $P$ is in $(x',y_n)$-form,
$P= \operatorname{OP}(p(x',y_n,\xi ))$. Since $p$ satisfies the
transmission condition, each term in the Taylor expansion in $y_n$
can be decomposed:
$$\aligned
p(x',y_n,\xi )&= \sum_{0\le l<l_0}\tfrac 1{l!}\partial _{y_n}^lp(x',0,\xi
)y_n^l+\tilde p_{l_0}(x',y_n,\xi )y_n^{l_0},\text{ with }\\
\partial _{y_n}^lp(x',0,\xi
)&= p_{(l)}'(x',\xi )+p_{(l)}''(x',\xi );
\endaligned\tag4.1$$
here  $\tilde p_{l_0}$ is a $\psi $do symbol of
order $\nu $, 
$p_{(l)}'$ is a differential operator symbol of order $\nu $ (void if $\nu
<0$) and  $p_{(l)}''$ is a $\psi $do symbol of order $\nu $ that is
$O(\ang{\xi _n}^{-1})$ in $\xi _n$ (cf\. (A.6)). We
take as usual 
$2k-\nu >n$.

\proclaim{Lemma 4.1} Each of the terms in the
decomposition
$$\multline
P_+G^{(k)}_\mu =P_+G_\mu ^{(k)}\theta\\
=\sum_{l<l_0}\tfrac 1{l!} [ \operatorname{OP}\bigl(p_{(l)}'(x',\xi
)\bigr)_+x_n^{l}G_\mu ^{(k)}\theta  +
\operatorname{OP}\bigl(p_{(l)}''(x',\xi  
)\bigr)_+x_n^{l}G_\mu ^{(k)}\theta  ]
\\+ \operatorname{OP}\bigl(\tilde p_{l_0}(x',y_n,\xi
)\bigr)_+x_n^{l_0}G_\mu ^{(k)}\theta \endmultline\tag4.2
$$
is trace class.

The trace norm of 
$ \operatorname{OP}\bigl(\tilde p_{l_0}(x',y_n,\xi
)\bigr)_+x_n^{l_0}G_\mu ^{(k)}\theta $ is
$O(\mu ^{n+1+\nu _+-2k-l_0})$, for each $l_0\in\Bbb N$ with
$l_0>n-2k+\nu _+$, $\nu _+=\max
\{\nu ,0\}$.

For each $l$, the trace of the operator $
\operatorname{OP}\bigl(p_{(l)}'(x',\xi 
)\bigr)_+x_n^{l}G_\mu ^{(k)}\theta $ has an asymptotic expansion in a series
without logarithms 
$$
\Tr  \operatorname{OP}\bigl(p_{(l)}'(x',\xi
)\bigr)_+x_n^{l}G_\mu ^{(k)}\theta \sim\sum_{r\in \Bbb N}c_r\mu ^{n-1+\nu
-2k-l-r}, \text{ for $\mu \to\infty $ in }\Gamma .\tag 4.3
$$

\endproclaim 
\demo{Proof} The decomposition (4.2) follows directly from (4.1),
 and
all the operators are trace class since they are of the form
$G'\theta =(\theta {G'}^*)^*$, where $G'$ and ${G'}^*$ are singular
Green operators of order $\nu -2k<-n$
and class 0.

For the second statement, we recall the rule
$$x_n^{l}\operatorname{OPG}(g)=\operatorname{OPG}(\overline D_{\xi
_n}^{l}g)\tag 4.4$$
(cf\. e.g\. [G96, Lemma 2.4.3]); thus $x_n^{l_0}G^{(k)}_\mu $ is
a singular Green operator that
is strongly polyhomogeneous of order $-2k-l_0$ and class 0. 
In  the terminology of [G96], it is of regularity $+\infty $,
whereas $\operatorname{OP}(\tilde p_{l_0})$, being of order $\nu $ and
$\mu $-independent, is of regularity $\nu $. The composed operator 
$\operatorname{OP}(\tilde p_{l_0})_+x_n^{l_0}G^{(k)}_\mu $ is then
an s.g.o\. of order $\nu -2k-l_0$, class 0 and regularity $\nu $, and
so is its adjoint, and they map $L_2(\rnp)$ into $H^{2k+l_0-\nu ,\mu
}(\crnp)$ with norm $O(\ang\mu ^{-\nu }+1)$ (cf\. [G96, Th\.s 2.5.11,
2.7.6]). Then the trace norm of $\operatorname{OP}(\tilde
p_{l_0})_+x_n^{l_0}G^{(k)}_\mu\theta  $ is $O(\ang\mu
^{n+1-2k-l_0+\nu _+})$.

For the third statement we note that since
$\operatorname{OP}(p_{(l)}')$ is a differential operator of order
$\nu $, and the
factor $x_n^l$ has the effect 
described in (4.4), the composite is a strongly
polyhomogeneous singular Green operator of order $\nu -2k-l$. Then
its trace has an expansion as in (3.6); this gives (4.3).
\qed  
\enddemo

It remains to study the terms$$
\tfrac1{l!} \operatorname{OP}\bigl(p_{(l)}''(x',\xi
)\bigr)_+x_n^lG^{(k)}_\mu \theta .\tag4.5  
$$
Fix $l$ and denote for simplicity $$
\operatorname{OP}(p_{(l)}''(x',\xi ))=P''=\operatorname{OP}(p''(x',\xi
)),\quad x_n^lG^{(k)}_\mu \theta =\operatorname{OPG}(g') .\tag4.6 
$$
Here we expand the symbol of $P''$ in a
convergent series as
in (A.6):
$$
p''(x',\xi ',\xi _n)=\sum_{m\in\Bbb Z}b_m(x',\xi
')\hat\varphi  ' _m(\xi
_n,\sigma ),\tag4.7
$$
where $b_m$ is rapidly decreasing in $S^{\nu +1}$.
Concerning $g'$, we note that $\overline D_{\xi _n}^lg^{(k)}$ has a
structure like that of $g^{(k)}$ except that the summation in
homogeneous terms starts with $J=l$; this must moreover be composed
with $\theta $, which again gives a symbol with the same structure. We can
take it in $y'$-form, considering
its expansion in strongly
homogeneous terms:$$
g'(y',\xi ,\eta _n,\mu )\sim \sum_{l\le J<\infty }g'_{-1-2k-J}(y',\xi
,\eta _n,\mu ) ;\tag4.8
$$
here the $J$'th term is of the form in (2.31) (with $x'$ replaced by
$y'$). As in the analysis around (3.1) it suffices to study $\Tr
P''_+\,\operatorname{OPG}(g'_{-1-2k-J})$ for each $J$. 

Now we proceed as in Section 3, writing$$\aligned
\Tr P''_+\,\operatorname{OPG}(g'_{-1-2k-J}) &=\Tr\operatorname{OPG}\bigl(p''(x',\xi )_+\circ_n 
g'_{-1-2k-J}(y',\xi ,\eta _n,\mu )\bigr)\text{ where}\\
p''\circ_n g'_{-1-2k-J}&=\sum_{m\in\Bbb
Z}g_{m,J}(x',y',\xi ,\eta _n,\mu )\text{ with}\\
g_{m,J}&=
b_m(x',\xi
')\hat\varphi '_m(\xi _n,\sigma )_+
\circ_n 
g'_{-1-2k-J}(y',\xi ,\eta _n,\mu ).
\endaligned\tag4.9
$$
The composition rule is:
$p(\xi
_n)_+\circ_n g'(\xi _n,\eta _n )=h^+_{\xi _n}[p(\xi
_n)g'(\xi _n,\eta _n )].$

Termwise, we have by (3.5),$$\aligned
\Tr\operatorname{OPG}(g_{m,J}(x',y',\xi ,\eta _n,\mu
))&=\Tr'\operatorname{OP}'( s_{m,J}(x',y',\xi ',\mu )),\text{
where}\\
 s_{m,J}(x',y',\xi ',\mu )&= \tr_ng_{m,J}(x',y',\xi ,\eta _n,\mu
).
\endaligned\tag 4.10
$$
So the crucial
step is to calculate the normal trace of the 
s.g.o\. arising from composing $\hat\varphi '_m(\xi _n,\sigma )
$ as a $\psi $do symbol with the fractions in (2.31).
Here
$$\aligned
s_{j,j',m}(x',\xi ',\mu )&=\tr_n \bigl( \hat\varphi '_m(\xi _n,\sigma )_+\circ_n\tfrac1{(\kappa ^++i\xi _n)^j(\kappa ^--i\eta
_n)^{j'}}\bigr)\\
&=\tr_n \bigl( h^+[\hat\varphi '_m(\xi  _n,\sigma
)
\tfrac1{(\kappa ^++i\xi _n)^j}]\tfrac 1{(\kappa ^--i\eta  
_n)^{j'}}\bigr)\\
&=\tfrac1{2\pi }\int  h^+[\hat\varphi '_m(\xi  _n,\sigma
)
\tfrac1{(\kappa ^++i\xi _n)^j}]\tfrac 1{(\kappa ^--i\xi   
_n)^{j'}}\,d\xi _n\\
&=\tfrac1{2\pi }\int  \hat\varphi '_m(\xi  _n,\sigma
)
\tfrac1{(\kappa ^++i\xi _n)^j}\tfrac 1{(\kappa ^--i\xi   
_n)^{j'}}\,d\xi _n,
\endaligned\tag4.11
$$
where we used that $h^-[\hat\varphi '_m(\xi  _n,\sigma
)
\frac1{(\kappa ^++i\xi _n)^j}]\frac 1{(\kappa ^--i\xi   
_n)^{j'}}$ is holomorphic on $\Bbb C_+$ and $O(\xi _n^{-2})$ for
$|\xi _n|\to\infty $.

\proclaim{Lemma 4.2} The symbol
$$
s_{j,j',m}=
\tfrac 1{2\pi }\int \frac{(\sigma -i\xi _n)^m}{(\sigma
+i\xi _n)^{m+1}}\frac
{1}{(\kappa ^++ i\xi
_n)^{j}(\kappa ^--
i\xi  
_n)^{j'}}\,d\xi _n
$$
satisfies for $m\ge 0$, $j$ and $j'\ge 1$:
$$s_{j,j',m}=\sum_{\Sb |m-m'|\le j''< j'\\  m'\ge
0\endSb}b_{jj'j''m'}(\kappa ^-)^{-j''}\frac{(\sigma -\kappa
^-)^{m'}}{(\sigma +\kappa ^-)^{m'+1}}(\kappa ^++\kappa ^-)^{-j-j'+1+j''},
\tag4.12$$
where the $b_{jj'j''m'}$ are universal constants that are
$O(m^{j'})$ for fixed $j,j'$.
This is a weakly polyhomogeneous symbol 
in $S^{-j-j',0}\cap S^{0,-j-j'}$.
There is a similar formula for $m\le 0$, with $j$ and $j'$
interchanged, $\kappa ^+$ and $\kappa ^-$
interchanged. \endproclaim 

\demo{Proof}
Let $m\ge 0$. Using (3.17) and (A.5) to reduce the
expressions $\xi _n^{j''}\partial _{\xi _n}^{j''}\hat\varphi '
_m$, we find:$$\aligned 
s_{j,j',m}
&=\tfrac 1{2\pi }\int \tfrac{(\sigma -i\xi _n)^m}{(\sigma
+i\xi _n)^{m+1}}\tfrac
{1}{(\kappa ^++ i\xi
_n)^{j}}\tfrac1{({j'}-1)!}(-i\partial _{\xi
_n})^{{j'}-1}\tfrac1{\kappa ^--i\xi  
_n}\,d\xi _n\\
&=\tfrac 1{2\pi }\int_{\Cal L_-}\sum_{0\le j''\le j'-1}\tbinom {j'-1}{j''}(i\partial _{\xi
_n})^{{j''}}\tfrac{(\sigma -i\xi _n)^m}{(\sigma
+i\xi _n)^{m+1}}(i\partial _{\xi
_n})^{{j'}-1-j''}\tfrac
{1}{(\kappa ^++ i\xi
_n)^{j}}\cdot\\
&\qquad\cdot \tfrac1{({j'}-1)!}\tfrac1{\kappa ^--i\xi  
_n}\,d\xi _n\\
&=
\sum_{0\le j''\le j'-1}\Bigl[(i\partial _{\xi
_n})^{{j''}}\tfrac{(\sigma -i\xi _n)^m}{(\sigma
+i\xi _n)^{m+1}}\Bigr]\tfrac
{c_{jj'j''}}{(\kappa ^++ i\xi
_n)^{j+j'-1-j''}}\Big|_{\xi _n=-i\kappa ^-}\\
&=
\sum_{0\le j''\le j'-1}{(-i\kappa ^-)^{-j''}}\Bigl[\xi _n^{j''}(i\partial _{\xi
_n})^{{j''}}\tfrac{(\sigma -i\xi _n)^m}{(\sigma
+i\xi _n)^{m+1}}\Bigr]\tfrac
{c_{jj'j''}}{(\kappa ^++ i\xi
_n)^{j+j'-1-j''}}\Big|_{\xi _n=-i\kappa ^-}\\
&=\sum_{\Sb |m-m'|\le j''< j'\\ m'\ge
0\endSb}b_{jj'j''m'}(\kappa ^-)^{-j''}\tfrac{(\sigma -\kappa
^-)^{m'}}{(\sigma +\kappa ^-)^{m'+1}}(\kappa ^++\kappa ^-)^{-j-j'+1+j''},
\endaligned$$
where the $b_{jj'j''m'}$ are universal constants that are $O(m^{j'})$
for fixed $j,j'$.
By the rules in [GS95], this is a wphg symbol in in $S^{-j-j',0}\cap
S^{0,-j-j'}$.  

There is a
similar calculation in case $m<0$ where the integral is instead
treated as a residue in $\Bbb C_+$. 
\qed\enddemo

We then find:

\proclaim{Proposition 4.3} The
 pseudodifferential symbols $
 s_{m,J}(x',y',\xi ',\mu )$
defined in {\rm (4.10)}
are in $S^{\nu -2k-J,0}(\Bbb R^{n-1}\times\Bbb R^{n-1},\Gamma )\cap
S^{\nu +1,-1 -2k-J}(\Bbb R^{n-1}\times\Bbb R^{n-1},\Gamma )$.

Moreover,
the sum over $m$ converges in all symbol norms, so that also the symbol$$
\multline
\tilde s_{(J)}(x',y',\xi ',\mu )=\sum_{m\in\Bbb Z}s_{m,J}(x',y',\xi
',\mu )\\
=\tr_n(p''(x',\xi )_+\circ_n g'_{-1-2k-J}(y',\xi ',\xi _n,\eta _n,\mu ))
\endmultline
\tag4.13
$$
is in 
$S^{\nu -2k-J,0}\cap
S^{\nu +1,-1-2k-J}$.

Finally,
$$\multline
\Tr\bigl(P''_+\,\operatorname{OPG}(
 g'_{-1-2k-J}(y',\xi ,\eta _n,\mu
))\bigr)
= \Tr'\operatorname{OP}'\bigl(\tilde s_{(J)}(x',y',\xi ',\mu )\bigr)\\
\sim \sum_{r=0}^\infty a^{4}_{J,r}\mu ^{n-1+\nu
-2k-J-r}+\sum_{r=0}^\infty (a^{5}_{J,r}\log\mu +a^{6}_{J,r})\mu ^{-1 -2k-J-r}.
\endmultline \tag4.14
$$

 \endproclaim 

\demo{Proof} 
The first statement follows from Lemma 4.2, when the
coefficients $b_m$ and the numerators in the expansions (2.31) are taken into
account. The summation in $m$ is handled as in the proof of Proposition 3.4,
since $b_m$ is rapidly decreasing in $m$ and the $m$-dependent
coefficients arising from the compositions in Lemma 4.2 are
polynomially bounded in $m$. The trace expansion also follows as in
Proposition 3.4.
\qed
\enddemo

Recall here that we have suppressed the indexation in $l$ referring
to the Taylor expansion of $p(x',y_n,\xi )$ in $y_n$, and that $J\ge l$ for
each $l$.

Again, as in the terms treated in Proposition 3.4, there is no
term of the form $c\mu ^{-2k}\log\mu $.

\subhead 5. Trace calculations for terms in $GQ^k_{\mu ,+}$ \endsubhead 

We shall finally treat the term $GQ^k_{\mu ,+}$, which, as we show,
gives a nontrivial contribution to the relevant log-coefficient.
Again, some manipulations will be performed before expansion of the
right hand factor in homogeneous terms as at the end of Section 2.

In view of the reduction to a coordinate patch, the term can be
written in the
form $GQ^{k}_{\mu ,+} $  
with $G=\operatorname{OPG}(g(x',\xi ,\eta _n))$,  $Q^{k}_\mu
=\operatorname{OP}(q^{\circ k}(y,\xi ,\mu )\theta (y_n))$, the
symbols vanishing 
for large $x '$ and $y'$, and $\theta \in C_{(0)}^\infty (\crp)$,
$\theta (y_n)=1$ near $y_n=0$. In principle, $q$ can be thought of
as determined as a 
parametrix of an elliptic symbol $p_{1}+\mu ^2$ defined for all
$y_n\ge 0$, with uniform estimates. The effect of the cutoff $\theta
$ can then be 
eliminated as follows:

\proclaim{Lemma 5.1} The operator $GQ^k_{\mu ,+}(1-\theta )$  is a
smoothing operator with a trace expansion$$
\sum_{r\in \Bbb N}c_r\mu ^{n-2k-r}\text{ for $\mu \to\infty $ in
}\Gamma .\tag5.1
$$
\endproclaim 

\demo{Proof} We use the old trick of ``nested commutators''. With
$\chi  _1(x_n)=1-\theta (x_n)$ on $\crp$, let $\chi _2$, $\chi _3$,
\dots be a sequence of $C^\infty $ functions supported in $\rp$ and
satisfying$$
\chi _i\chi _{i+1}=\chi _i,\text{ for all }i,\tag5.2
$$ 
in other words, $\chi _{i+1}$ is 1 on $\operatorname{supp}\chi _i$
(but is still 0 on a neighborhood of 0). Denote $Q^k_{\mu }=Q_{(1)}$.
Now$$
Q_{(1)}  \chi _1=Q_{(1)}  \chi _1\chi _2=\chi _1Q_{(1)} \chi
_2+[Q_{(1)} ,\chi _1]\chi _2.
$$
Here $G\chi _1Q_{(1) ,+}\chi _2=G\chi _1(Q_{(1) }\chi _2)_+$, where
$G\chi _1$ is a $\mu $-independent operator of order $-\infty $, so
that the composed 
operator is 
the restriction to $\rnp$ of an operator 
in $S^{-\infty ,-2k}(\Bbb R^{n}\times\Bbb R^n,\Gamma )$;
its trace has an expansion as in (5.1) by [GS95, Th\. 2.1] (since $G$
maps into functions supported in a fixed compact set).

Let $Q_{(1)}$ have the symbol $q_{(1)}(y,\xi ,\mu )$; then the
commutator $Q_{(2)}=[Q_{(1)},\chi _1]$ has the  symbol$$ 
q_{(1)}(y,\xi ,\mu )(\chi _1(y_n)-\chi _1(x_n))=
q_{(1)}(y,\xi ,\mu )(y_n-x_n)\tilde \chi _1(x_n,y_n),
$$
where $\tilde\chi _1(x_n,y_n)$ is a smooth, bounded function with
bounded derivatives. The
$\psi $do with this symbol acts in the same way as the $\psi $do with
the symbol $D_{\xi _n}q_{(1)}(y,\xi ,\mu )\chi _1(x_n,y_n)$, which is
strongly polyhomogeneous of degree $-2k-1$. It also has a symbol
$q_{(2)}$ in $y$-form, which is  strongly polyhomogeneous of degree
$-2k-1$.
So $Q_{(2)}$ is one degree better than $Q_{(1)}$. We can continue the
argument, writing$$
Q_{(i)}  \chi _i=Q_{(i)}  \chi _i\chi _{i+1}=\chi _iQ_{(i)} \chi
_{i+1}+[Q_{(i)} ,\chi _i]\chi _{i+1}, \quad Q_{(i+1)}= [Q_{(i)} ,\chi _i],
\tag5.3$$
in the general step, finding that $G\chi _iQ_{(i)}\chi _{i+1}$
contributes with a trace expansion like (5.1) but starting with the
power $\mu ^{n-2k-i}$, and $Q_{(i+1)}$ is strongly polyhomogeneous of
degree $-2k-i-1$. The latter gives an operator with trace norm
$O(\mu ^{n+\nu _+-2k-i})$ for large enough $i$, so since we can let
$i\to\infty $ we get by superposition an expansion (5.1) for the
original operator.\qed
\enddemo 

\subsubhead{5.a The $y_n$-independent case}\endsubsubhead

For the rest of the analysis we first consider the case where $P_1$,
when in $y$-form, is independent of $y_n$ near
$\{y_n=0\}$, in the local coordinates we are using, so that we can
assume that $Q^k_\mu =\operatorname{OP}(q^{\circ k}(y',\xi ,\mu ))$.
As usual, it suffices to consider the 
contributions from the homogeneous 
terms in $q^{\circ k}(y',\xi ,\mu )\sim \sum_{J\in \Bbb
N}q_{-2k-J}^{\circ k}(y',\xi ,\mu )$. Moreover, we expand the symbol
of $G$ as in (3.11). Then
$$\aligned
\Tr G\,\operatorname{OP}(q^{\circ k}_{-2k-J}(y',\xi ,\mu ))_+&=\Tr\operatorname{OPG}(g(x',\xi ,\eta _n)\circ_n 
q^{\circ k}_{-2k-J}(y',\xi ,\mu )),\text{ where}\\
g\circ_n q^{\circ k}_{-2k-J}&=\sum_{l,m\in\Bbb
N}g'_{l,m,J}(x',y',\xi ,\eta _n,\mu )\text{ with}\\
g'_{l,m,J}&=
2\sigma c_{lm}(x',\xi
')\hat\varphi '_l(\xi _n,\sigma )
\bar{\hat\varphi }'_m(\eta _n,\sigma )\circ_n 
q^{\circ k}_{-2k-J}(y',\xi ,\mu )_+.
\endaligned\tag5.4
$$The composition rule is: $g(\xi _n,\eta _n)\circ_n q(\xi
_n)_+=h^{-}_{\eta _n}[g(\xi _n,\eta _n) q(\eta 
_n)]$.

Termwise, we have by (3.5),$$\aligned
\Tr\operatorname{OPG}(g'_{l,m,J}(x',y',\xi ,\eta _n,\mu
))&=\Tr'\operatorname{OP}'(\tilde s_{l,m,J}(x',y',\xi ',\mu )),\text{
where}\\
\tilde s_{l,m,J}(x',y',\xi ',\mu )&= \tr_ng'_{l,m,J}(x',y',\xi ,\eta _n,\mu
),
\endaligned\tag 5.5
$$
and the crucial
step is to calculate the normal trace of the 
s.g.o\. arising from composing $\hat\varphi '_l(\xi _n,\sigma )
\bar{\hat\varphi }'_m(\eta _n,\sigma )$ with the fractions in (2.10)
(decomposed as in (2.15)).
In view of (3.9),$$
\multline
\tr_n(
2\sigma c_{lm}\hat\varphi '_{l}(\xi _n,\sigma )\bar{\hat\varphi
}'_{m}(\eta  _n,\sigma )\circ_n q^{\circ k}_{-2k-J})\\
=\tfrac1{2\pi }\int 
2\sigma c_{lm}\hat\varphi '_{l}(\xi _n,\sigma )\bar{\hat\varphi
}'_{m}(\xi  _n,\sigma )q^{\circ
k}_{-2k-J}(y',\xi ',\xi _n,\mu )\,d\xi _n.
\endmultline
\tag5.6
$$
It is seen from the
decomposition of $q^{\circ
k}_{-2k-J}(y',\xi ',\xi _n,\mu )$ in Proposition 2.1 that the main point is
to calculate the expressions
$$
s^{\pm}_{j,l,m}(y',\xi ',\mu )=\tfrac 1{2\pi }\int \frac 
{(\sigma - i\xi
_n)^{l}}{(\sigma + i\xi
_n)^{l+1}}\frac{(\sigma +
i\xi  
_n)^{m}}{(\sigma -
i\xi  
_n)^{m+1}}\frac1{(\kappa
^{\pm}\pm i\xi _n)^j}\,d\xi _n.
\tag5.7
$$

\proclaim{Lemma 5.2} One has for $l,m\in\Bbb Z$, $j\ge 1$:$$\alignedat2
&\text{For $m<l$, }&&s^+_{j,l,m}=0,\quad
s^-_{j,l,m}=\sum_{ |l-m-1-m'|<j,\,  m'\ge 
0}(\kappa^-) ^{1-j}b'_{jm'}\frac{(\sigma -\kappa ^-)^{m'}}{(\sigma +\kappa ^-
)^{m'+2}},\\
&\text{For $m=l$, }&&s^+_{j,l,l}=\frac1{(\kappa ^++\sigma )^j 2\sigma
},\quad s^-_{j,l,l}=\frac1{(\kappa ^-+\sigma )^j 2\sigma },\\
&\text{For $m>l$, }&&s^+_{j,l,m}=\sum_{ |m-l-1-m'|<j,\, m'\ge
0}(\kappa ^+) ^{1-j}b_{jm'}\frac{(\sigma -\kappa
^+)^{m'}}{(\sigma +\kappa ^+ 
)^{m'+2}},\quad s^-_{j,l,m}=0,
\endalignedat\tag5.8$$
with $b_{jm'}$ and $b'_{jm'}$ being $O(l^jm^j)$ for fixed $j$.

The resulting symbols when $m\ne l$ are in $S^{-j-1,0}\cap
S^{0,-j-1}$,
whereas the symbols for $m=l$
are in  $S^{-j-1,0}\cap S^{-1,-j}$.

\endproclaim 

\demo{Proof}
Consider $s^+_{j,l,m}$. The integrand
is rational with poles in $\Bbb C\setminus\Bbb R$ and is $O(\xi
_n^{-3})$ for $|\xi _n|\to\infty $, so the integral can be deformed to 
an integral over a large closed curve in $\Bbb C_+$ or $\Bbb C_-$.
$$
\text{For $m<l$, }s^+_{j,l,m}=\tfrac 1{2\pi }\int \frac{(\sigma - 
i\xi  
_n)^{l-m-1}}{(\sigma +
i\xi  
_n)^{l-m+1}}\frac1{(\kappa
^{+}+ i\xi _n)^j}
\,d\xi _n=0,\tag5.9
$$
since the integrand is holomorphic on ${\Bbb C}_-$.

When $m=l$, we use that there is just the simple pole $-i\sigma $ in
$\Bbb C_-$: 
$$
\text{For $m=l$, } s^+_{j,l,l}=\tfrac 1{2\pi }\int \frac1{(\sigma +
i\xi  
_n)(\sigma -i\xi _n)(\kappa
^{+}+i\xi _n)^j}
\,d\xi _n=\frac1{2\sigma (\kappa ^++\sigma )^j },\tag5.10
$$
considering the integral as a residue at $-i\sigma $.
This function is the product of a symbol in $S^{-j,0}\cap S^{0,-j}$
and a $\mu $-independent symbol
in $S^{-1,0}$, hence it
lies in  $S^{-j-1,0}\cap S^{-1,-j}$.

When $m>l$, we use instead that the only pole in $\Bbb C_+$ is
$i\kappa ^+$. Writing $m-l-1=k\ge 0$, we find:
$$\aligned
s^+_{j,l,m}&=\tfrac 1{2\pi }\int \frac{(\sigma + 
i\xi  
_n)^{m-l-1}}{(\sigma -
i\xi  
_n)^{m-l+1}}\frac1{(\kappa
^{+}+ i\xi _n)^j}
\,d\xi _n\\
&=\tfrac 1{2\pi }\int_{\Cal L_+}
\bar{\hat\varphi}' _{k}(\xi _n,\sigma )\bar{\hat\varphi}' _{0}(\xi
_n,\sigma )
\frac1{(\kappa
^{+}+ i\xi _n)^j}
\,d\xi _n\\
&=\tfrac 1{2\pi }\int_{\Cal L_+}
\bar{\hat\varphi}' _{k}(\xi _n,\sigma )\bar{\hat\varphi}' _{0}(\xi
_n,\sigma )
\tfrac1{(j-1)!}(i\partial _{\xi _n})^{j-1}\frac1{\kappa ^++i\xi 
_n}
\,d\xi _n\\
&=\tfrac1{2\pi
}\int_{\Cal L_+}\Bigl[\tfrac1{(j-1)!}(-i\partial _{\xi
_n})^{j-1}\bigl(\bar{\hat\varphi }'_{k}(\xi _n,\sigma
)\bar{\hat\varphi }'_{0
}(\xi _n,\sigma )\bigr)\Bigr]\frac1{\kappa ^++i\xi
_n}
\,d\xi _n\\
&=\tfrac{(-i)^{j-1}}{(j-1)!}
\partial _{\xi _n}^{j-1}\bigl[\bar{\hat\varphi }'_{k}(\xi _n,\sigma
)\bar{\hat\varphi }'_{0}(\xi _n,\sigma
)\bigr]\Big|_{\xi _n=i\kappa ^+}\\
&=\tfrac{(-\kappa ^+)^{1-j}}{(j-1)!}
\xi _n^{j-1}\partial _{\xi _n}^{j-1}\bigl[\bar{\hat\varphi }'_{k}(\xi
_n,\sigma )\bar{\hat\varphi }'_{0}(\xi _n,\sigma
)\bigr]\Big|_{\xi _n=i\kappa ^+}\\
&=\sum_{\max\{k-j+1,0\}\le m'\le
k+j-1}(\kappa^+) ^{1-j}b_{jm'}\frac{(\sigma -\kappa ^+)^{m'}}{(\sigma +\kappa ^+
)^{m'+2}},\\
\endaligned\tag5.11$$
very similarly to the calculations in (3.18);
the $b_{jm'}$ are $O(l^jm^j)$. 
This is a symbol in $S^{-j-1,0}\cap S^{0,-j-1}$.

The study
of $s^-_{j,l,m}$ is similar, except that the roles of the upper and
lower half-planes, as well as the roles of $m$ and $l$, are
interchanged.
\qed
\enddemo

Because of the difference between the cases, we sum over $l=m$ and
$l\ne m$ separately: 

\proclaim{Proposition 5.3} 
The pseudodifferential symbols $
 \tilde s_{l,m,J}(x',y',\xi ',\mu )$
defined in {\rm (5.5)} are in
$S^{\nu -2k-J,0}\cap
S^{\nu +1,-1 -2k-J}$ when $
l\ne m$ and in
$S^{\nu -2k-J,0}\cap
S^{\nu ,-2k-J}$ when $l=m$.

Moreover,
the summations over $l\ne m$ and over $l=m$ converge in the respective
symbol norms, so that 
$$
\aligned
\tr_n(g(x',\xi ',\xi _n,\eta _n)&\circ_n q^{\circ k}_{-2k-J}(y',\xi
',x_n,\mu )) =\tilde s_{(J)}(x',y',\xi ',\mu )\\
&=\tilde s^{\prime}_{(J)}(x',y',\xi ',\mu )+\tilde s^{\prime\prime}_{(J)}(x',y',\xi ',\mu )
,\text{ with }\\
\tilde s^{\prime}_{(J)}(x',y',\xi ',\mu )&=\sum_{l\in\Bbb
N}s_{l,l,J}(x',y',\xi 
',\mu )\in S^{\nu -2k-J,0}\cap
S^{\nu ,-2k-J},\\
\tilde s^{\prime\prime}_{(J)}(x',y',\xi ',\mu )&=\sum_{
l,m\in\Bbb N,\, l\ne m}\tilde s_{l,m,J}(x',y',\xi
',\mu )\in S^{\nu -2k-J,0}\cap
S^{\nu +1,-1-2k-J}.\\
\endaligned\tag5.12
$$
 
Finally,
$$\aligned
\Tr\bigl(G\,\operatorname{OP}(
q^{\circ k}_{-2k-J}&(y',\xi ,\mu
))\bigr)
= 
\Tr'\operatorname{OP}'\bigl(\tilde s'_{(J)}\bigr)+\Tr'\operatorname{OP}'\bigl(\tilde s''_{(J)}\bigr)\\
&\sim \sum_{l=0}^\infty a^{7}_{J,l}\mu ^{n-1+\nu
-2k-J-l}+\sum_{l=0}^\infty (a^{8}_{J,l}\log\mu +a^{9}_{J,l})\mu ^{ -2k-J-l},\\
\text{with }\Tr'\operatorname{OP}'\bigl(\tilde s^{\prime}_{(J)}\bigr)
&\sim \sum_{l=0}^\infty a^{7\prime}_{J,l}\mu ^{n-1+\nu
-2k-J-l}+\sum_{l=0}^\infty (a^{8\prime}_{J,l}\log\mu
+a^{9\prime}_{J,l})\mu ^{ -2k-J-l},\\ 
\Tr'\operatorname{OP}'\bigl(\tilde s^{\prime\prime}_{(J)}\bigr)
&\sim \sum_{l=0}^\infty a^{7\prime\prime}_{J,l}\mu ^{n-1+\nu
-2k-J-l}+\sum_{l=0}^\infty (a^{8\prime\prime}_{J,l}\log\mu +a^{9\prime\prime}_{J,l})\mu ^{-1 -2k-J-l}.\endaligned\tag5.13
$$
In particular, the only contribution to a term $c\mu ^{-2k}\log\mu $ comes from
$\tilde s^{\prime}_{(0)}$.
\endproclaim 

\demo{Proof} 
The first statement follows from Lemma 5.2, when the
coefficients $c_{lm}(x',\xi ')$ and the numerators in the expansions
(2.15) are taken into 
account. The summations in $l$ and $m$ are handled as in the proof of
Proposition 3.4, using the
rapid decrease of the $c_{lm}$ in $l$ and $m$ and the polynomial
bounds on the $(l,m)$-dependent
coefficients arising from the compositions in Lemma 5.2. Then the
trace expansions follow as in 
Proposition 3.4.
\qed
\enddemo 

We shall now study the
log-contribution from $\tilde s_{(0)}^{\prime}$.

\proclaim{Proposition 5.4} We have (with $\sigma =\sigma (\xi ')=[\xi
']$):
$$\multline
\tr_n(
2\sigma c_{ll}\hat\varphi '_{l}(\xi _n,\sigma )\bar{\hat\varphi
}'_{l}(\eta  _n,\sigma )\circ_n q^{\circ k}_{-2k})\\
=c_{ll}(x',\xi ')(q^{\circ k,+}_{-2k}(y',\xi ',-i\sigma
,\mu )+q^{\circ k,-}_{-2k}(y',\xi ',i\sigma ,\mu )),
\endmultline \tag5.14
$$
lying in $S^{\nu -2k,0}\cap S^{\nu ,-2k}$. Summation over
$l$ gives the symbol$$
\tilde s_{(0)}^{\prime}(x',y',\xi ',\mu )=(\tr_n g)(x',\xi ')\cdot
(q^{\circ k,+}_{-2k}(y',\xi ',-i\sigma
,\mu )+q^{\circ k,-}_{-2k}(y',\xi ',i\sigma,\mu )) .\tag 5.15
$$
\endproclaim
\demo{Proof} 
In
view of (5.6)--(5.8) and Proposition 2.1 (in particular (2.14)), the
left hand side of 
(5.14) equals  
$$
\aligned
\tfrac1{2\pi }\int  2\sigma c_{ll}
&\hat\varphi
'_l(\xi _n,\sigma )\bar{\hat\varphi }'_l(\xi _n,\sigma )\sum_{1\le j\le
k}\bigl(\frac{a^{+}_{k,j}} 
{a^k(\kappa ^++i\xi _n )^{j}}+\frac{a^{-}_{k,j}}
{a^k(\kappa ^--i\xi _n )^{j}}\bigr)\,d\xi _n\\
&=\frac{2\sigma c_{ll}}{2\sigma }\sum_{1\le j\le k}\Bigl(\frac{a^{+}_{k,j}}
{a^k(\kappa ^++ \sigma )^{j}}+\frac{a^{-}_{k,j}}
{a^k(\kappa ^-+\sigma )^{j}}\Bigr)\\
&=c_{ll} (q^{\circ k,+}_{-2k}(y',\xi ',-i\sigma
,\mu )+q^{\circ k,-}_{-2k}(y',\xi ',i\sigma ,\mu ))
.\endaligned\tag5.16
$$ This shows (5.14), and (5.15) follows since
$$
\aligned
\tr_n g(x',\xi ',\xi _n,\eta _n)&=
\tr_n \sum_{l,m}2\sigma c_{lm}(x',\xi ')\hat\varphi '_l(\xi_n,\sigma
)\bar{\hat\varphi} '_m(\eta 
_n,\sigma )\\
&=\sum_{l,m} 2\sigma c_{lm}(x',\xi ')(\varphi '_l,\varphi
'_m)=\sum_{l}c_{ll}(x',\xi '),\endaligned \tag5.17
$$
where we use that $\frac1{2\pi }(\hat\varphi '_l,\hat\varphi '_m)=(\varphi '_l,\varphi '_m)=\frac1{2\sigma }\delta _{lm}$.
\qed
\enddemo 

The contribution to the log-term in (5.13) is found as follows:
Denote $$
q^{\circ
k,+}_{-2k}(x',\xi ',-i\sigma ,\mu  )+q^{\circ
k,-}_{-2k}(x',\xi ',i\sigma ,\mu 
)=\alpha (x',\xi ',\mu ).\tag5.18
$$According to [GS95, pf\. of Th\. 2.1], we have to calculate (on the
diagonal, where $x'=y'$)
$$\multline
f(x',\mu )=\int_{\Bbb R^{n-1}}(\tr_ng)_{1-n}(x',\xi ')\alpha (x',\xi ',\mu )\,\d\xi'\\=\int_{|\xi '|\ge |\mu
|}(\tr_ng)_{1-n}\alpha \,\d\xi'+\int_{|\xi '|\le
1}(\tr_ng)_{1-n}\alpha \,\d\xi'+\int_{1\le |\xi '|\le |\mu 
|}(\tr_ng)_{1-n}\alpha \,\d\xi' \endmultline\tag5.19
$$ at each $x'$ (denoting $(2\pi )^{1-n}d\xi '=\d\xi '$); it is integrated afterwards in $x'$. 
The integrand is homogeneous of degree $1-n-2k$ in
$(\xi ',\mu )$ for 
$|\xi '|\ge 1$, so  the first integral equals a
constant times $\mu ^{-2k}$ (for $|\mu |\ge 1$). It is seen as in
[GS95] that the second integral (over 
$|\xi '|\le 1$) contributes to the pure powers $c_l\mu ^{-l}$. It is
the third integral that produces a logarithm. To find the
coefficient, we need to study $\alpha (x',\xi ',\mu )$ more closely. Since it
is in $S^{0,-2k}$, it has a ``Taylor expansion at $\infty $'' by
[GS95, Th\. 1.12]. More
precisely this means that
 for $z=\frac1\mu $ lying in $\Gamma $ (2.7),
 $z^{-2k}\alpha (x',\xi ',\Zfrac)$ has a
Taylor expansion  at $z= 0$:$$
\alpha (x',\xi ',\Zfrac)\sim z ^{2k}(\alpha _0(x',\xi ')+\alpha _1(x',\xi ')z+\dots+
\alpha _r(x',\xi ')z^{r}+\dots)
\text{ for $z \to0 $ in
}\Gamma ,
\tag5.20$$
which we also write$$
\multline
\alpha (x',\xi ',\mu )\sim \alpha _0(x',\xi ')\mu ^{-2k}+\alpha _1(x',\xi ')\mu 
^{-2k-1}+\dots+
\alpha _r(x',\xi ')\mu ^{-2k-r}+\dots\\
\text{ for $\mu \to\infty $ in
}\Gamma 
.\endmultline\tag5.21
$$
Here $\alpha _r\in S^r$, and the remainder after the $r$'th
term is in $S^{r+1, -2k-r-1}$. 

We shall show below in Lemma 5.5 that $\alpha _0(x',\xi ')=1$. 
Then 
$$\multline
\int_{1\le |\xi '|\le |\mu
|}(\tr_ng)_{1-n}\alpha \,\d\xi' \\
\sim \int_{1\le |\xi '|\le |\mu
|}(\tr_ng)_{1-n}(x',\xi ')(\mu ^{-2k}+\alpha _1(x',\xi ')\mu
^{-2k-1}+\dots)\,\d\xi' . \endmultline\tag5.22
$$Here, since $(\tr_ng)_{1-n}(x',t\eta
')=t^{1-n}(\tr_ng)_{1-n}(x',\eta ')$ for 
$|\eta '|=1$,$$\aligned
\int_{1\le |\xi '|\le |\mu
|}(\tr_ng)_{1-n}&(x',\xi ')\mu ^{-2k}\,\d\xi' \\
&=\tfrac1{(2\pi )^{n-1}}\mu
^{-2k}\int_{t=1}^{t=\mu }t^{-1}\,dt\int_{|\eta '|=1}(\tr_ng)_{1-n}(x',\eta
')\,d\sigma (\eta ')\\
&=  \mu
^{-2k}\log \mu \tfrac1{(2\pi )^{n-1}}\int_{|\eta '|=1}(\tr_ng)_{1-n}(x',\eta
')\,d\sigma (\eta ') \endaligned\tag5.23 
$$
(where $d\sigma (\eta
')$ indicates the surface measure). The other terms give
contributions of the form $c_r(x')\mu ^{-2k-r}$, $r>0$, and the
remainder after $r$ terms is $O(\mu ^{-M})$ where $M\to\infty $ for
$r\to\infty $, as shown in [GS95]. We can conclude that the
coefficient of 
$\mu ^{-2k}\log\mu $ in (5.13) is
$$
a^{8\prime}_{0,0}=\tfrac1{(2\pi )^{n-1}}\int_{\Bbb R^{n-1}}\int_{|\xi
'|=1}(\tr_ng)_{1-n}(x',\xi 
')\,d\sigma (\xi ')dx' .\tag5.24
$$

\proclaim{Lemma 5.5} 
In the expansion
{\rm (5.21)} of the function $\alpha (x',\xi ',\mu )=q^{\circ
k,+}_{-2k}(x',\xi ',-i\sigma ,\mu  )+q^{\circ k,-}_{-2k}(x',\xi ',i\sigma ,\mu 
)$, the coefficient $\alpha _0(x',\xi ')$ of $\mu ^{-2k}$ equals $1$.
\endproclaim

\demo{Proof} 
We have  that
$$
p_{1,2}(x',0,\xi )+\mu ^2=a(x')\xi _n^2+b(x',\xi ')\xi _n+c(x',\xi ')+\mu
^2,\tag5.25 
$$
where $\operatorname{Re}a(x')>0$, $b(x',\xi ')$ is linear in $\xi '$,
$c(x',\xi ')$ is quadratic in $\xi '$. The fact that the roots $\pm i\kappa
^\pm(x',\xi ',\mu )$ with respect to $\xi _n$
lie in $\Bbb C_{\pm}$, respectively, when $\mu \in \Gamma $, follows
from the strong ellipticity; the
precise formula for $\kappa ^{\pm}$ is $$ 
\kappa ^{\pm}=(\mu ^2/a+c/a-(b/2a)^2)^{\frac12}\mp ib/2a.\tag5.26
$$
The functions $\kappa ^\pm(x',\xi ',\mu )$ are strongly homogeneous of
degree 1 in $(\xi ',\mu )$, so they have
``Taylor expansions for $\mu \to\infty $'' in $\Gamma $ (by [GS95, Th\. 1.12]):$$
\kappa ^\pm(x',\xi ',\mu )\sim a(x')^{-\frac12}\mu +\varrho^{\pm} _1(x',\xi ')\mu ^0+\dots+\varrho
^\pm_r(x',\xi ')\mu ^{1-r}+\dots\tag5.27
$$
with $\varrho ^\pm_r\in S^r$; the first coefficient is
$a^{-\frac12}$ in view of (5.26).
If we add a $\mu $-independent term $t\sigma (\xi ')$, $t\in[0,1]$,
the functions 
$\kappa ^{\pm}+t\sigma $ 
likewise have expansions 
$$
\kappa ^\pm(x',\xi ',\mu )+t \sigma (\xi ')\sim a(x')^{-\frac12}\mu +(\varrho
^\pm_1(x',\xi ')+t\sigma (\xi '))\mu ^0+\dots+\varrho ^\pm  
_r(x',\xi ')\mu ^{1-r}+\dots,\tag5.28
$$
{\it with the same first term}. Since 
$\operatorname{Re}\kappa ^\pm>0$, they have inverses
$(\kappa ^\pm+t\sigma )^{-1}$
with
$$
\multline
(\kappa ^\pm(x',\xi ',\mu )+t\sigma (\xi '))^{-1}
 \sim 
a(x')^{\frac12}\mu ^{-1}+\tilde\varrho ^\pm_1(x',\xi ',t)\mu
^{-2}\\+\dots+\tilde\varrho ^\pm  
_r(x',\xi ',t)\mu ^{-1-r}+\dots,
\text{ for $\mu \to\infty $ in }\Gamma ,\endmultline\tag5.29
$$
 where the first coefficient is
necessarily the inverse of the
preceding first coefficient, again {\it independent of $t$}.
 Now consider the decomposition (2.14) for $m=k$:
$$\aligned
q^{\circ k}_{-2k}&=\frac1{(p_{1,2}(x',\xi )+\mu ^2)^{k}}=\frac1{a^k(\kappa
^++i\xi
_n)^{k}(\kappa ^--i\xi
_n)^{k}}=q^{\circ k,+}_{-2k}+q^{\circ k,-}_{-2k},\\
q^{\circ k,\pm}_{-2k}&=\sum_{1\le j\le k}\frac{a^\pm_{k,j}}{a^k(\kappa
^\pm \pm i\xi 
_n)^{j}}.\endaligned\tag5.30
$$
Here the numerators $a_{k,j}^\pm(x',\xi ',\mu )$ 
are strongly homogeneous in $(\xi ',\mu )$ of degree $j-2k$, hence have
expansions (by [GS95, Th\. 1.12])
$$\multline
a^{\pm}_{k,j}(x',\xi ',\mu )
\sim a^{\pm}_{kj0}(x',\xi ')\mu ^{j-2k}+
a^{\pm}_{kj1}(x',\xi ')\mu ^{j-2k-1}\\+\dots +a^{\pm}_{kjr}(x',\xi ')\mu
^{j-2k-r}+\dots,
\text{ for }\mu \to\infty \text{ in }\Gamma .\endmultline \tag5.31
$$
Inserting these expansions as well as those from (5.29) in the formulas
(5.30) for $q^{\circ k,\pm}_{-2k}$,
we find the following expansions:
$$\multline
q^{\circ k,\pm}_{-2k}(x',\xi ',\mp it\sigma (\xi '),\mu )=\sum_{1\le j\le
k}\frac{a^\pm_{k,j}(x',\xi ',\mu )}{a(x')^k(\kappa
^\pm(x',\xi ',\mu )+t\sigma (\xi '))^{j}}\\
\sim 
\alpha _0^\pm(x',\xi ')\mu ^{-2k}+\alpha _1^{\pm}(x',\xi ',t)\mu
^{-2k-1}+\dots+\alpha _r^\pm 
(x',\xi ',t)\mu ^{-2k-r}+\dots,\endmultline\tag5.32
$$
where the first coefficient $\alpha _0^\pm(x',\xi ')$ is {\it independent
of }$t\in [0,1]$. 
For $t=1$ it follows from (5.18), (5.21) and (5.32) that$$
\alpha _0(x',\xi ')=\alpha ^+_0(x',\xi ')+\alpha ^-_0(x',\xi ').\tag5.33
$$
On the other hand, we have for
$t=0$ that
$$
q^{\circ k,+}_{-2k}(x',\xi ',0,\mu )+q^{\circ k,-}_{-2k}(x',\xi
',0,\mu )=q^{\circ k}_{-2k}(x',\xi ',0,\mu ),  
$$
where 
$$\multline
q^{\circ k}_{-2k}(x',\xi ',0,\mu )=(p_{12}(x',\xi ',0)+\mu ^2)^{-k}
\\
\sim \mu ^{-2k}+\tilde \alpha _1(x',\xi ',0)\mu ^{-2k-1}+\dots+
\tilde \alpha _r(x',\xi ',0)\mu ^{-2k-r}+\dots\text{ for }\mu \to\infty
\text{ in }\Gamma ,
\endmultline
$$
with coefficient 1 for the first term.
Hence in view of (5.32),$$
\alpha ^+_0(x',\xi ')+\alpha ^-_0(x',\xi ')=1.\tag5.34
$$
(5.33) and (5.34) together show the lemma.
\qed
\enddemo

\smallskip

\subsubhead{5.b The $y_n$-dependent case}\endsubsubhead

Finally, we consider the case where $p_1(y,\xi )$ and $q(y,\xi ,\mu
)$ do depend on $y_n$ near $\{y_n=0\}$. Denote $p_1(y',0,\xi
)=p_1^0(y',\xi )$ and $q(y',0,\xi
,\mu )=q^0(y',\xi ,\mu )$. Here $(y',y_n)$ is used in a compact
subset of $\Bbb R^n$, but we can assume that the symbols have extensions to
$\Bbb R^n$ defining operators $P_1$, $P^0_1$, $Q_\mu $, $Q^0_\mu $ so
that $Q_\mu =(P_1+\mu ^2)^{-1}$, $Q^0_\mu =(P^0_1+\mu ^2)^{-1}$. Then
$$\aligned
Q_\mu -Q^0_\mu&=Q_\mu (P^0_1+\mu ^2)Q^0_\mu  -Q_\mu (P_1+\mu
^2)Q^0_\mu =Q_\mu (P^0_1-P_1)Q^0_\mu,\\  
(Q_\mu)^k -(Q^0_\mu)^k
&=(Q_\mu)^k -(Q_\mu )^{k-1}Q^0_\mu+(Q_\mu
)^{k-1}Q^0_\mu-\dots\\
&\qquad-Q_\mu (Q^0_\mu)^{k-1} +Q_\mu (Q^0_\mu)^{k-1} 
-(Q^0_\mu )^k\\
&=
(Q_\mu )^k(P^0_1-P_1)Q^0_\mu+ (Q_\mu )^{k-1}(P^0_1-P_1)(Q^0_\mu)^2+\dots\\
&\qquad +Q_\mu (P^0_1-P_1)(Q^0_\mu)^{k}. 
\endaligned\tag5.35$$ 
In view of the
formulas (2.10) for the symbols of $Q_\mu $ and $Q_\mu ^0$, we find
that the symbol $\tilde q(y,\xi ,\mu )$ of $(Q_\mu)^k -(Q^0_\mu)^k$
has the structure of a series $\sum_{J\in\Bbb N}\tilde
q_{-2k-J}(y,\xi ,\mu )$ of  homogeneous terms of the form$$
\tilde q_{-2k-J}(y,\xi ,\mu )=\sum_{J/2+k+1\le m+m'\le 2J+k+1}\frac{r_{k,J,m,m'}(y,\xi
)}{(p_{1,2}(y,\xi )+\mu ^2)^{m}(p^0_{1,2}(y,\xi )+\mu ^2)^{m'}}\tag 5.36
$$
with $r_{k,J,m,m'}(y,\xi )$ denoting a homogeneous polynomial in $\xi $ of
degree $2(m+m')-2k-J$.

As usual, it suffices to consider the contribution from each $\tilde
q_{-2k-J}$. 
This must be reduced to situations where the composed symbols can be
calculated. We can write$$
p^0_{1,2}(y',\xi )-p_{1,2}(y,\xi )=y_n\tilde p(y,\xi ),\tag 5.37
$$
where $\tilde p$ is a polynomial in $\xi $ of degree 2. Using the
formula $1/{(1-b)}=\sum_{0\le j<N}b^j+b^N/(1-b)$, valid for $b\ne 1$,
we write, for $y_n$ so small that $|y_n\tilde p/(p^0_{1,2}+\mu ^{2})|\le
c<1$,$$\aligned
\frac1{p_{1,2}+\mu ^{2}}&=\frac1{p^0_{1,2}+\mu ^{2}} \Bigl(
\frac{p_{1,2}+\mu ^{2}}{p^0_{1,2}+\mu ^{2}}\Bigr)^{-1}=
\frac1{p^0_{1,2}+\mu ^{2}} \Bigl(
1-\frac{y_n\tilde p}{p^0_{1,2}+\mu ^{2}}\Bigr)^{-1}\\
&=\sum_{0\le j<N}\frac{y_n^j\tilde p^j}{(p^0_{1,2}+\mu ^{2})^{j+1}}
+\frac{y_n^N\tilde p^N}{(p^0_{1,2}+\mu ^{2})^N(p_{1,2}+\mu ^{2})}.
\endaligned\tag5.38$$
For the last term in (5.38) we observe that when $N$ is even,$$
\frac{\tilde p^N}{(p^0_{1,2}+\mu ^{2})^N(p_{1,2}+\mu ^{2})}=\tilde
p^{N/2}\Bigl(\frac{\tilde p}{(p^0_{1,2}+\mu
^{2})^2}\Bigr)^{N/2}\frac1{p_{1,2}+\mu ^{2}}\equiv \tilde
p^{N/2}\tilde p_{(N)},  
$$
where $\tilde p^{N/2}\in S^N$, whereas 
$\tilde p_{(N)}$    belongs to $S^{-N-2,0}\cap 
S^{0,-N-2}$.

For the study of the composed operator
$G\operatorname{OP}\bigl(r(y,\xi )y_n^N \tilde p^{N/2}(y,\xi )
\tilde p_{(N)}(y,\xi ,\mu )\bigr)_+$, we can ``disentangle'' the $\mu
$-independent factor
$r(y,\xi )y_n^N\tilde p^{N/2}$ from the $\mu $-dependent factor
$\tilde p_{(N)}$ as follows:
We observe that the ordinary product of two pseudodifferential symbols
$p$ and $p'$ given in $y$-form may be
expressed as a series of compositions of derived symbols:
$$
p(y,\xi )p'(y,\xi ) \sim \sum_{\alpha
\in\Bbb N^{n}}\tfrac1{\alpha !}\partial _{y}^\alpha  p(y,\xi
)\circ D_{\xi }^\alpha p'(y,\xi ).\tag 5.39
$$
This ``disentanglement formula'' (which can be deduced from the
usual formula for composition of symbols in $y$-form by successive
comparison of terms) shows
how a $\psi$do, given in the form $\operatorname{OP}(pp')$, may be written as
a series of truly composed $\psi $do's.

Applying  (5.39), the first term gives the composite $G\operatorname{OP}(r(y,\xi
)y_n^N\tilde p^{N/2})_+\operatorname{OP}(\tilde p_{(N)})_+$, where the
order of $G\operatorname{OP}(r(y,\xi
)y_n^N\tilde p^{N/2})_+$ equals $\nu $ + the order of $r$,
thanks to the fact that $y_n^N$ reduces the order of a s.g.o\. by
$N$, compensating for the factor $\tilde p^{N/2}$. The factor
$\operatorname{OP}(\tilde p_{(N)})$ here gives lower order and better $O(\mu
^{-M})$ estimates, the larger $N$ is taken. There is also a leftover
term (of the type $GL(Q',Q'')$), that behaves similarly. In the other
contributions resulting from the use of (5.39), the degree in $y_n$
goes down, but this is compensated for by a lower degree of the right
hand
factor. The remainder after the terms up to $|\alpha |=l$ in the
expansion has better 
$\mu $-estimates, the larger $l$ is taken.

In reality, the last term in (5.38) enters in (5.36) in powers, composed with
terms having only powers of $p^0_{1,2}+\mu ^2$ in the denominator,
but the pattern of the resulting operators is similar to what we have
just
described.

This leaves us with contributions of the form $r(y,\xi
)(p^0_{1,2}(y',\xi )+\mu ^2)^{-m}$. When $r(y,\xi )$ is disentangled 
from $(p^0_{1,2}(y',\xi )+\mu ^2)^{-m}$ by (5.39) as above and the resulting
operators are composed with $G$, we arrive at the type of
operators we have dealt with in the first part of this section. They
do have trace expansions like the first one in (5.13), when of order
$\nu -2k-J$. For 
$J>0$ these do not contribute to the coefficient of
$\mu ^{-2k}\log\mu $. 

For $J=0$, we moreover have to show that the resulting trace does not
contribute to the $\mu ^{-2k}\log \mu $-term. Here we have to deal
with the principal symbol of (5.35), equal to$$
\sum_{i=1}^k\frac{y_n\tilde p}{(p^0_{1,2}+\mu ^2)^i(p_{1,2}+\mu
^2)^{k+1-i}}
=
\sum_{i=1}^k\frac{\sum_{l,m=1}^ny_na_{lm}\xi _l\xi _m}{(p^0_{1,2}+\mu
^2)^{k+1}}\bigl(\sum_{0\le j<N}\frac{y_n^j\tilde p^j}{(p^0_{1,2}+\mu ^2)^j}\bigl)^{k-i}
+  \text{ rem.;}
$$ 
we have applied
(5.38) to the factors $1/(p_{1,2}+\mu ^2)^{k+1-i}$ and inserted
$\tilde p=\sum_{l,m}a_{lm}(y)\xi _l\xi _m$. For each $l,m$, we 
can move a factor $y_na_{lm}(y)\xi _l$ over to $G$ in the same way as above
without changing the order and class, so we are left with
a finite sum of compositions of s.g.o.'s of order $\nu $ and class $0$
with $\psi 
$do's with symbols$$
\frac{\xi _m}{(p^0_{1,2}+\mu
^2)^{k+1}}\Bigl(\frac{y_n\tilde p}{p^0_{1,2}+\mu ^2}\Bigl)^{j(k-i)}.
$$
We decompose these rational functions of $\xi _n$ using (2.14) and
the rules$$
\frac {i\xi _n}{p^0_{1,2}+\mu ^2}=-\frac{\frac{\kappa ^+}{\kappa ^++\kappa
^-}}{\kappa 
^++i\xi _n}+\frac{\frac{\kappa ^-}{\kappa ^++\kappa ^-}}{\kappa 
^--i\xi _n},\quad
\frac {\xi _n^2}{p^0_{1,2}+\mu ^2}=1-\frac{\frac{(\kappa ^+)^2}{\kappa ^++\kappa
^-}}{\kappa 
^++i\xi _n}-\frac{\frac{(\kappa ^-)^2}{\kappa ^++\kappa ^-}}{\kappa 
^--i\xi _n},
$$
and find that the resulting simple fractions are $O(\mu ^{-2k-1})$; in the
calculation of $\tr_n$ of the
compositions with the s.g.o.'s (as in Proposition 5.3) this leads to
weakly polyhomogeneous $\psi $do symbols 
in $S^{\nu +1,-2k-1}$ on $\Bbb R^{n-1}$. They do not contribute to
$\mu ^{-2k}\log \mu $.

We have shown:

\proclaim{Proposition 5.6} Each term in the expansion
$G((Q_\mu
)^k-(Q^0_\mu )^k)_+\sim \sum_{J}G\,\operatorname{OP}(\tilde q_{-2k-J})_+$ 
has a trace expansion as in {\rm (5.13)}, with vanishing coefficient of $\mu
^{-2k}\log\mu $.
\endproclaim

\subhead 6. Synthesis  \endsubhead

We now have all the ingredients for the proof of Theorem 1.1, in the
following detailed formulation:

\proclaim{Theorem 6.1}
Let $P$ be a classical $\psi $do of order $\nu \in\Bbb Z$ defined on
a neighborhood $\widetilde X$ of $X$ and 
having the transmission property at $X'$, and let $G$ be a singular
Green operator of order $\nu$ and class $0$ on $X$.
The operator
$P$ acts in a $C^\infty $ vector bundle $\widetilde E$ over
$\widetilde X$, and
 $P_+$ and $G$ act in $E=\widetilde E|_{X}$.

Let $P_1$ be
a strongly 
elliptic second-order differential operator in $\widetilde E$ (with scalar
principal symbol on a neighborhood of $X'$) and let
$P_{1,\operatorname{D}}$ be its Dirichlet realization on $X$; by adding a
constant morphism we can assume that $P_1$ and
$P_{1,\operatorname{D}}$ have positive lower bound, so that the
resolvents $(P_1-\lambda )^{-1}$ and
$(P_{1,\operatorname{D}}-\lambda )^{-1}$ exist for $-\lambda $ in a
region $W$ {\rm (2.1)}.
We can assume that 
$\widetilde X$ is compact, without boundary.

Let $k>(n+\nu )/2$. Then there are full asymptotic trace expansions {\rm
(1.10)}, with 
the coefficients $\tilde c_j, \tilde c_l',\tilde c''_l$ 
proportional to the coefficients $c_j,  c_l', c''_l$ by
universal constants. 
In {\rm
(1.10 I)} $\lambda \to\infty $ on the rays in $-W$, in {\rm (1.10
II)} $t\to 0+$, and {\rm (1.10 III)} holds in the sense that the
left hand side minus the terms with $l<N$, is holomorphic
for $\operatorname{Re}s>-N/2$, for any $N$.
 Here $\tilde c'_0= c'_0$ and equals
$$
\tfrac1{2(2\pi )^n}\int_{S^*X}\tr_E p_{-n}(x,\xi )\,d\sigma
(\xi )dx+\tfrac 1{2(2\pi )^{n-1}}\int_{S^*X'}\tr_E(\tr_ng)_{1-n}(x',\xi
')\,d\sigma (\xi ')dx'.
\tag6.1$$ 
Thus the {\rm [FGLS96]} residue of $P_++G$ recalled in {\rm (1.5)}
satisfies $$
\operatorname{res}(P_++G)=2c'_0=\operatorname{ord}P_1\cdot \operatorname{Res}_{s=0}\Tr((P_++G)P_{1,\operatorname{D}}^{-s}).\tag6.2
$$

\endproclaim 

\demo{Proof} 
Since $P$ can be taken to be 0 outside a neighborhood of $X$, and
$P_1$ has to be strongly elliptic (the symbol homotopic to $|\xi
|^2I$), there are no topological obstructions to 
obtaining a compact $\widetilde X$.

Consider the five terms in (2.5). 
The four s.g.o\. terms have in Section 2 been partitioned into finite
sums of terms 
supported in coordinate patches, with errors that are either $O(\mu
^{-N})$, any $N$, or contribute trace expansions as in (2.8). This
allowed us to work in $\crnp$, where we formulated the results for
the model case $\dim E=1$. We
have shown at the end of Section 2 that the contributions from the
remainder after $J_0$ terms in the polyhomogeneous expansions of the
symbols of  $Q^k_{\mu ,+}$, $G^-(Q^k_\mu )$ and $G^{(k)}_\mu $
give operators whose trace norms are $O(\mu ^{-N})$ with
$N\to\infty $ for $J_0\to \infty $. Similar results hold for other
remainders encountered along the way, such as the last term in (4.2)
containing the Taylor remainder
(when $l_0\to\infty $), the
error in cut-offs treated in Lemma 5.1, 
and the errors that stem from expansions (5.38) (when $N\to\infty $).
There remain, up to a given accuracy $O(\mu ^{-M})$ in the desired
trace expansion, a finite number (depending on $M$) of contributions,
that have been described in  Propositions 3.4, 4.3, 5.3 and 5.6.
It is important to observe here that the expansions stemming from the
precise terms 
$p'_{(l)}$ and $p''_{(l)}$ in the
Taylor expansion (4.1)
give trace expansions beginning with the monomial $\mu ^{n-1+\nu -2k-l}$
and the log-monomial $\mu ^{-1-2k-l}\log\mu $, with powers going to
$-\infty $ for $l\to\infty $, so that any specific term gets only finitely
many contributions when we let $l_0\to\infty $.

When $\dim E>1$, the above considerations pertain to each element in
the matrix compositions.

Adding all the terms, we find an expansion:
$$
\Tr [ -L(P,Q^k_\mu )+P_+G^{(k)}_\mu +GQ^k_{\mu ,+}+GG^{(k)}_\mu ]
\sim \sum_{j\ge 0} a_j\mu ^{
n+\nu -2k-j}+\sum_{l\ge 0}  (a'_{l}\log\mu +a''_{l})\mu ^{ -2k-l}.
$$
Here it is only
$GQ^k_{\mu ,+}$ that contributes to $a'_0$, by expressions (5.24) for
the diagonal elements in 
each localized term.
Insertion of $\mu ^2=-\lambda $ gives
$$\multline
\Tr [ -L(P,Q^k_\mu )+P_+G^{(k)}_\mu +GQ^k_{\mu ,+}+GG^{(k)}_\mu ]\\
\sim
\sum_{j\ge 0}  a_j(-\lambda ) ^{\frac {n+\nu -j}2-k}+\sum_{l\ge 0}
(\tfrac12a'_{l}\log(-\lambda ) +a''_{l})(-\lambda 
) ^{-\frac l2 -k}.\endmultline\tag6.3
$$
Now, as
observed already in [W84] and [FGLS96], (5.24) has a meaning independent of
the choice of local coordinates, so the coefficient of $(-\lambda
)^{-k}\log(-\lambda )$ satisfies
$$
\tfrac12 a'_0=\tfrac 1{2(2\pi )^{n-1}}\int_{S^*X'}\tr_E(\tr_ng)_{1-n}(x',\xi
')\,d\sigma (\xi ')dx'.\tag6.4
$$

It remains to account for the first expression in (2.5). Here 
[GS95, Th\. 2.1] implies an expansion of the kernel on the diagonal:
$$
\Cal K_{PQ_\mu ^k}(x,x)\sim
\sum_{j\ge 0} b_j(x)\mu ^{n+\nu
-2k-j}+\sum_{l\ge0} (b'_{l}(x)\log\mu +b''_{l}(x))\mu ^{ -2k-l};\tag6.5
$$
and an investigation of the proof (as in Section 5 above) shows that
since $$
q_{-2k}=\frac1{(p_{1,2}+\mu ^2)^{k}}=\mu ^{-2k}
(1+p_{1,2}\mu ^{-2})^{-k}\sim \mu 
^{-2k}(1+\varrho _2\mu ^{-2}+\varrho _4\mu ^{-4}+\dots)\tag6.6
$$ 
{\it with first coefficient} 1, the contribution from each diagonal
element to the coefficient of $\mu
^{-2k}\log\mu $ in local coordinates has the form $\frac1{(2\pi
)^n}\int_{|\xi |=1}p_{-n}(x,\xi )\,d\sigma (\xi )$. 
Summation of the pieces and integration over $X$ gives:
$$
\Tr (PQ_\mu ^k)_+\sim
\sum_{j\ge 0}  b_{j,+}\mu ^{n+\nu 
-2k-j}+\sum_{l\ge0} (  b'_{l,+}\log\mu + b''_{l,+})\mu ^{
-2k-l}.\tag6.7 
$$
Again we can use
the invariance shown in [W84] and [FGLS96] to justify writing$$
 b'_{0,+}=\tfrac 1{(2\pi )^{n}}\int_{S^*X}\tr_E p_{-n}(x,\xi
)\,d\sigma (\xi )dx.\tag6.8
$$
This is turned into an expansion in powers of $-\lambda $ as above.
Collecting the contributions, we obtain (1.10), with $ c'_0$ described
by (6.1).

The other expansions (1.10 II) and (1.10 III) are derived from this
by the transition principles worked out e.g\. in [GS96],
based on the 
formulas (1.8). One checks from [GS96, Cor\. 2.10] that $\tilde c'_0=c'_0$. 
\qed

\enddemo 

\subhead Appendix. Laguerre expansions \endsubhead

Here we recall the basic properties of the
Laguerre function systems used in [G96], and introduce a slight
modification. 
We first recall the definitions (these and the listed properties are
worked out in detail in [G96, Sect\. 2.2]):$$
\aligned
\widehat\varphi_k
(\xi _n,\sigma)&=(2\sigma
)^{\frac12}\frac{(\sigma-i\xi _n)^k}{(\sigma+i\xi _n)^{k+1}},\quad
k\in\Bbb Z,\\ 
\varphi_k(x_n,\sigma)&=
 (2\sigma)^{\frac 12}H(x_n)(\sigma-\partial_{x_n})^k(x_n^ke^{-x_n\sigma})/k!
\text{ for }k\ge 0,\\
\varphi_k(x_n,\sigma)&=\varphi_{-k-1}(-x_n,\sigma)\text{ for }k<0.
\endaligned\tag A.1
$$
($H(x_n)$ is the Heaviside function $1_{x_n>0}$.)
Here $\sigma $ can be any positive number, and the systems $\{\varphi
_k\}_{k\in\Bbb Z}$ resp\. $\{\varphi
_k\}_{k\in\Bbb N}$ are orthonormal bases of $L_2(\Bbb R)$ resp\.
$L_2(\Bbb R_+)$. 
The $\varphi_k$ with $k\ge 0$ are
the
eigenfunctions of the (unconventional) Laguerre operator
 $$
\Cal
L_{\sigma,+}=\sigma^{-1}(\sigma+\partial_x)x(\sigma-\partial_x)=
-\sigma^{-1}\partial_x x\partial_x+\sigma x+1\tag A.2
 $$
 in $L_2(\Bbb R_+)$, with simple eigenvalues $2(k+1)$, and the
$\varphi_k$ with $k<0$ are similarly the eigenfunctions for $\Cal
L_{\sigma,-}$ 
defined by the same expression on $\Bbb R_-$.

When $u\in L_2(\Bbb R_+)$ is expanded in the
Laguerre system 
$(\varphi_k)_{k\in\Bbb N}$, by
 $$
u(x_n)=\sum_{k\in\Bbb N} b_k\varphi_k({x_n},\sigma),\tag A.3
 $$
then $u\in\Cal S(\overline{\Bbb R}_+)=r^+\Cal S(\Bbb R)$ if and only if
$(b_k)_{k\in\Bbb N}$ is rapidly decreasing.

The norming factor $\sqrt{2\sigma }$ in (A.1) can be inconvenient in
some calculations, so let us also introduce a notation for the
functions without it:$$
\aligned
\widehat\varphi'_k
(\xi _n,\sigma)&=\frac{(\sigma-i\xi _n)^k}{(\sigma+i\xi
_n)^{k+1}}=(2\sigma )^{-\frac12}\widehat\varphi_k(\xi _n,\sigma ),\\ 
\varphi'_k(x_n,\sigma)&=(2\sigma )^{-\frac12}\varphi _k(x_n,\sigma).
\endaligned\tag A.4
$$They satisfy$$
\alignedat2
& \partial_{\xi _n}\hat\varphi'_k({\xi _n},\sigma) =
\tfrac{-i}{2\sigma}(k\hat\varphi'_{k-1}+(2k+1)\hat\varphi'_k+
(k+1)\hat\varphi'_{k+1}), && \text{ for $k\in\Bbb Z$},\\
& i{\xi _n}\hat\varphi'_k({\xi _n},\sigma) =
-\sigma\hat\varphi'_k+2\sigma\tsize\sum_{0\le j< k}
(-1)^{k-1-j}\hat\varphi'_j+(-1)^k, && 
\text{ for $k\ge 0$},\\
&\partial _{\xi _n}(\xi _n\hat\varphi '_k(\xi _n,\sigma))=-\tfrac
k2\hat\varphi '_{k-1}+\tfrac12 \hat\varphi '_k
+\tfrac{k+1}2\hat\varphi '
_{k+1}, && \text{ for $k\in\Bbb Z$},\\
&\partial _\sigma \hat\varphi '_k(\xi _n,\sigma )=\tfrac1{2\sigma
}(k\hat\varphi '_{k-1}-\hat\varphi '_k-(k+1)\hat\varphi '_{k+1}), && \text{ for $k\in\Bbb Z$}.
\endalignedat\tag A.5
$$
The formulas show how simple manipulations with Laguerre functions lead to
expressions involving neighboring Laguerre functions.

When $p(x',0,\xi )$ is of order $d\in \Bbb Z$ and satisfies the
transmission condition, then$$
p(x',0,\xi )
=\sum_{0\le l\le d}s_{l}(x',\xi ')\xi
_n^l+\sum_{k\in\Bbb Z}b_k(x',\xi ')\hat\varphi '_k(\xi _n,\sigma );
\tag A.6
$$ 
the $s_l$ are polynomials in $\xi '$ of degree $d-l$ and the $b_k$
form a sequence that is rapidly decreasing in $S^{d+1}$ for $|k|\to\infty $.

Finally, we note two useful formulas:$$
(i\partial _{\xi _n})^k\frac1{\sigma +i\xi _n}=
\frac { k!}{(\sigma +i\xi _n)^{k+1}},\quad
\Cal F^{-1}\frac 1{(\sigma +i\xi _n)^{k+1}}=H(x_n)\frac{x_n^k}{k!}e^{-\sigma x_n}.\tag A.7
$$

\enddocument